\documentclass[final]{siamart251104}
\newif\ifacm
\newif\ifarxiv
\arxivtrue

\usepackage[utf8]{inputenc}
\usepackage[english]{babel}

\ifacm
\usepackage{balance} % For balanced columns on the last page
\else
\usepackage{amsmath}
\usepackage{amsfonts}
\usepackage{amssymb}
\usepackage{amsopn}
\fi

\ifacm
\else
\fi

\ifacm
\else
\usepackage{booktabs}
\fi
\usepackage{tabularx}

\usepackage{color}
\usepackage{bm}
\usepackage{lineno}
\usepackage{url}
\usepackage{hyperref}
\usepackage{cleveref}
\usepackage[boxruled,algo2e]{algorithm2e}
\usepackage{setspace}
\usepackage{xcolor}
\usepackage{wrapfig}
\usepackage{graphicx}
\usepackage{pstool}
\usepackage{multirow}
\usepackage{setspace}
\usepackage{enumitem}
\usepackage{multirow}
\usepackage{subcaption}
\usepackage{comment}
\usepackage{xspace} %usefull for defining command such as \fftw, and avoid trailing space, see https://tex.stackexchange.com/a/290504
\usepackage{tablefootnote} % Footnote in tables

\graphicspath{{./figures/}}

\DeclareCaptionFormat{custom}
{%
    \textbf{#1#2}\textnormal{\small #3}
}
\usepackage{pgfplots}
\usepgfplotslibrary{external} 
\pgfplotsset{compat=newest}

\usepackage{tikz}% pour diagrammes
\usetikzlibrary{shapes,snakes}
\usetikzlibrary{patterns}
\usetikzlibrary{shapes.misc}
\usetikzlibrary{arrows}
\usepgfplotslibrary{groupplots}
\usepgfplotslibrary{dateplot}
\usetikzlibrary{backgrounds}
\usetikzlibrary{external}
\usetikzlibrary{decorations.pathreplacing,calligraphy}

\newcommand{\code}[1]{\texttt{#1}}

\newcommand{\flups}[1]{\code{FLUPS}#1}
\newcommand{\murphy}[1]{\code{murphy}#1}
\newcommand{\pforest}[1]{\code{p4est}#1}
\newcommand{\hn}{\hat{\mathbf{n}}}

\newcommand{\normi}[1]{\vert\!\vert #1 \vert\!\vert_{\infty}}

\newcolumntype{Y}{>{\centering\arraybackslash}X}
\newcolumntype{Z}[1]{>{\hsize=#1\arraybackslash}X}

\newcommand{\be}{\begin{equation}}
\newcommand{\eec}{\;\;,\end{equation}}
\newcommand{\ee}{ \end{equation}}
\newcommand{\eed}{\;\;.\end{equation}}
\newcommand{\bse}{\begin{subequations}}
\newcommand{\ese}{\end{subequations}}

\newcommand{\Fref}[1]{\Cref{#1}}

\newcommand{\fref}[1]{\cref{#1}}
\newcommand{\eqqref}[1]{\cref{#1}}
\newcommand{\algref}[1]{\cref{#1}}
\newcommand{\sref}[1]{\cref{#1}}

\definecolor{darkGreen}{RGB}{0,160,0}
\definecolor{darkBlue}{RGB}{51,51,204}
\definecolor{darkRed}{RGB}{160,0,0}
\definecolor{orange}{RGB}{255,120,0}
\definecolor{myGreen}{rgb}{0.0706, 0.5333, 0.1921}

\newcommand{\revone}[1]{}
\newcommand{\revthree}[1]{}
\newcommand{\revtwo}[1]{}
\newcommand{\revall}[1]{}
\newcommand{\typo}[1]{}

\definecolor{airforceblue}{rgb}{0.36, 0.54, 0.66}
\definecolor{carolinablue}{rgb}{0.6, 0.73, 0.89}
\definecolor{moonstoneblue}{rgb}{0.45, 0.66, 0.76}
\definecolor{amber}{rgb}{1.0, 0.75, 0.0}
\definecolor{cadmiumorange}{rgb}{0.93, 0.53, 0.18}
\definecolor{persianorange}{rgb}{0.85, 0.56, 0.35}
\definecolor{cadetgrey}{rgb}{0.57, 0.64, 0.69}
\definecolor{darkgray}{rgb}{0.66, 0.66, 0.66}
\definecolor{asparagus}{rgb}{0.53, 0.66, 0.42}
\definecolor{cambridgeblue}{rgb}{0.64, 0.76, 0.68}
\definecolor{olivine}{rgb}{0.6, 0.73, 0.45}
\definecolor{carnelian}{rgb}{0.7, 0.11, 0.11}
\definecolor{cardinal}{rgb}{0.77, 0.12, 0.23}
\definecolor{firebrick}{rgb}{0.7, 0.13, 0.13}

\definecolor{matlabBlue}{rgb}{0, 0.4470, 0.7410}
\definecolor{matlabOrange}{rgb}{0.8500, 0.3250, 0.0980}
\definecolor{matlabYellow}{rgb}{0.9290, 0.6940, 0.1250}
\definecolor{matlabPurple}{rgb}{0.4940, 0.1840, 0.5560}
\definecolor{matlabGreen}{rgb}{0.4660, 0.6740, 0.1880}
\definecolor{matlabRed}{rgb} {0.7765, 0.1569, 0.0902}
\definecolor{matlabBlack}{rgb}{0, 0, 0}

\definecolor{keyBlue}{rgb}{0.0392, 0.2274, 0.4235}
\definecolor{keyBlueTrans}{rgb}{0.4274, 0.4980, 0.5764}
\definecolor{keyGreen}{rgb}{0.05490, 0.38039, 0.003921}
\definecolor{keyGreenTrans}{rgb}{0.7019, 0.8313, 0.7176}
\definecolor{keyOrange}{rgb}{0.98039, 0.4980, 0.03529}
\definecolor{keyGray}{rgb}{0.50196, 0.50196, 0.50196}
\definecolor{keyGreen2}{rgb}{0.0627 , 0.4196 , 0.3882}%{16,107,99}
\definecolor{keyOrange2}{rgb}{0.9843 , 0.2941 , 0.2431}%{255,99,78}

\definecolor{pyBlue}{rgb}{0.11372,0.42352,0.67059}
\definecolor{pyBlueLight}{rgb}{0.6196,0.79215,0.88235}
\definecolor{pyOrange}{rgb}{1.0,0.4549,0.0627}
\definecolor{pyOrangeLight}{rgb}{0.9921568,0.643137,0.376470}
\definecolor{pyGreen}{rgb}{0.1529,0.5882,0.1529}
\definecolor{pyRed}{rgb}{0.81569,0.13725,0.14117}
\definecolor{pyPurple}{rgb}{0.53725,0.36078,0.7098}
\definecolor{myCacaDoie}{rgb}{0.580392156862745,0.43,0.27}
\definecolor{pyBrown}{rgb}{0.549019607843137,0.337254901960784,0.294117647058824}
\definecolor{pyRose}{rgb}{0.890196078431372,0.466666666666667,0.76078431372549}
\definecolor{pyOrange2}{rgb}{0.9921,0.5098,0.2078}
\definecolor{pyBlue2}{rgb}{0.3765,0.6431,0.8157}
\definecolor{pyBar1}{rgb}{0.1647,0.2745,0.2745}
\definecolor{pyBar2}{rgb}{0.0000,0.4588,0.4588}
\definecolor{pyBar3}{rgb}{0.4863,0.7804,0.9998}
\definecolor{pyBar4}{rgb}{0.4549,0.4549,0.4549}

\tikzset{cross/.style={cross out, draw, minimum size=2*(#1-\pgflinewidth), inner sep=0pt, outer sep=0pt}}

\definecolor{cadmiumorange}{rgb}{0.93, 0.53, 0.18}
\definecolor{asparagus}{rgb}{0.53, 0.66, 0.42}
\definecolor{cerulean}{rgb}{0.0, 0.48, 0.65}
\definecolor{brickred}{rgb}{0.8, 0.25, 0.33}
\definecolor{ferngreen}{rgb}{0.31, 0.47, 0.26}
\hypersetup{
    colorlinks=true,
    linkcolor=cerulean,
    filecolor=magenta,      
    urlcolor=brickred,
    citecolor=ferngreen
    }
\newcommand{\picfont}{\sffamily}
\pgfplotsset{
	reverse legend,
	legend style={anchor=north west,at={(0.01,0.99)},font=\tiny},
	legend image code/.code={
		\draw[very thick,
			 mark repeat=2,
			 mark phase=2,
			 dash phase=0pt,
			 ] plot coordinates {(0pt,0pt) (10pt,0pt)};%
	}
	}
\tikzset{font=\picfont{}}
\newcommand{\per}[1]{\,\!/\,\!}
\newcommand{\uspace}[1]{\,\!}
\definecolor{col_rma_single}{RGB}{214,39,40}
\definecolor{col_single}{RGB}{255,127,14}
\definecolor{col_rma_multi}{RGB}{227,119,194}
\definecolor{col_part}{RGB}{31,119,180}
\definecolor{col_multi}{RGB}{99,121,57}
\definecolor{col_rma_active}{RGB}{189,158,57}
\definecolor{col_rma_single_active}{RGB}{176,176,176}
\definecolor{col_stream}{RGB}{188,189,34}

\definecolor{col_model_black}{RGB}{176,176,176}
\definecolor{col_model_manyvci}{RGB}{255,127,14}
\definecolor{col_model_manyvcitwo}{RGB}{99,121,57}
\definecolor{col_model_onevci}{RGB}{31,119,180}

\title{A scalable high-order multigrid-FFT Poisson solver for unbounded domains on adaptive multiresolution grids}

\author{
    Gilles Poncelet \thanks{Institute of Mechanics, Materials and Civil Engineering, Université catholique de Louvain, Belgium (corresponding author: \email{gilles.poncelet@uclouvain.be}).}
    \and Jonathan Lambrechts \footnotemark[1]
    \and Thomas Gillis \thanks{Department of Mechanical Engineering, Massachusetts Institute of Technology, Cambridge, MA, United States. \textit{Present address:} NVIDIA, Boulder, CO, United States.}
    \and Philippe Chatelain \footnotemark[1]
}

\begin{document}

\maketitle
\begin{abstract}
    Multigrid solvers are among the most efficient methods for solving the Poisson equation, which is ubiquitous in computational physics. For example, in the context of incompressible flows, it is typically the costliest operation. The present document expounds upon the implementation of a flexible multigrid solver that is capable of handling any type of boundary conditions within \murphy, a multiresolution framework for solving partial differential equations (PDEs) on collocated adaptive grids. The utilization of a Fourier-based direct solver facilitates the attainment of flexibility and enhanced performance by accommodating any combination of unbounded and semi-unbounded boundary conditions. The employment of high-order compact stencils contributes to the reduction of communication demands while concurrently enhancing the accuracy of the system. The resulting solver is validated against analytical solutions for periodic and unbounded domains. In conclusion, the solver has been demonstrated to demonstrate scalability to 16,384 cores within the context of leading European high-performance computing infrastructures.
\end{abstract}


%!TEX root = umg_1_intro.tex
%!TEX encoding = UTF-8 Unicode

\section{Introduction}
\label{section:intro}
The solution of Poisson problems is ubiquitous in computational physics as it concerns problems ranging from electromagnetism to fluid dynamics. When considering time-dependent applications, one has to solve this problem many times per simulation, making it one of the costliest operations. Developing efficient Poisson solvers is therefore a critical challenge and is tackled by many scientists nowadays. In this work, we consider the Poisson equation in a 3D domain:
\begin{equation}
    \nabla^2u= f\,. \label{eq:poisson_problem}
\end{equation}
In order to form a closed problem, this needs to be complemented with a set of boundary conditions. Those can include unbounded boundary conditions, also referred to as free-space boundary conditions, which suppose that the solution goes to zero at infinity. This is typical of problems involving a potential field induced by sources located in a compact region of space.

A traditional approach to solving the Poisson equation with unbounded directions has relied on the following approximation: the solution decay at an infinite distance is replaced by the enforcement of homogeneous Dirichlet conditions at a finite distance, taken as large as possible. 
However, this approach tends to be either computationally expensive or inaccurate when attempting to capture the solution decay by using a very large computational domain.
Alternative techniques make use of the Green's function convolution to accurately solve unbounded problems on compact domains \cite{hockney_computer_2021}. Typically, the choice of fast convolution algorithm is made between the fast Fourier transform (FFT) \cite{caprace_flups_2021, hejlesen_high_2013, chatelain_fourier-based_2010, gabbard_lattice_2023} or fast multipole method \cite{liska2014, hou_adaptive_2024, gillman_fast_2010}. Multigrid methods can also be combined with an unbounded solver to enable unbounded boundary conditions \cite{teunissen_geometric_2019}.

A most challenging situation occurs when the Poisson solution is sought for sources that exhibit a wide range of scales; this is classically the case in Computational Fluid Dynamics (CFD): external aerodynamics problems involve both very fine but crucial scales within boundary layers and much larger coherent structures shed by the device, i.e. its wake. As a consequence, computational requirements vary greatly across the domain if one aims to achieve a given error. Furthermore, given the nature of these problems, they are often very expensive to solve and require efficient parallel implementations that can scale up to very large partitions of High-Performance Computing infrastructures \cite{ibeid_fft_2020}.

Among the many families of Poisson solvers, this work will focus on the variants that rely on Cartesian grids, and especially those involving multiresolution adaptation. Those are indeed integral to the targeted CFD methods, namely block-based multi-resolution discretizations \cite{gillis_murphy:2022}, finite differences \cite{ji_fourth_2025} or hybrid Lagrangian-Eulerian methods such as the Vortex Particle-Mesh method \cite{chatelain_vortex_2008, balty_2025}. 
If one considers the Cartesian grid requirement, potential solvers include the already mentioned fast Fourier transform(FFT)-based techniques, fast multipole methods and geometric multigrid methods.

FFT-based techniques are arguably the most efficient at solving the Poisson equation on uniform grids while accommodating a flexible choice of boundary conditions, also including unbounded directions, see e.g. \cite{caprace_flups_2021}. There are however strictly limited to uniform discretizations.
Complex flow physics, with their broad range of scales, then render this approach inefficient and costly as one has to enforce the finest resolution required throughout the domain. On the other hand, multipole and multigrid methods directly address such a complexity and exhibit better performance when sources are strongly localized \cite{gholami_fft_2016}.
Furthermore, one might argue that multipole methods are fundamentally oblivious to structured discretizations. Therefore, as this contribution targets an implementation within an adaptive multiresolution Cartesian framework, we will focus on the multigrid approach for the remainder of this discussion.
%A multigrid approach was chosen over a multipole method due to pre-existing data structure fitting the former better. Coupling such a solver with an adaptive multiresolution framework then enables large-scale simulation of multiscale physical phenomena. 

Multigrid methods stem from the observation that the convergence rate of smoothers such as Gauss-Seidel or Jacobi depend on the wavenumber content of the solution, with  typically, a more effective smoothing of high wavenumbers.
The core principle of the multigrid methodology then lies in the efficient handling of low wavenumber modes by a smoother on a coarser mesh \cite{trottenberg_multigrid_2000}.
Thence, the iterative application over a hierarchy of successively coarser meshes can efficiently smoothen all the modes of the error.
This produces a solution technique, which, when applied to a discretized elliptic PDE such as \eqqref{eq:poisson_problem}, converges to a solution in an optimal $\mathcal{O}(N)$ fashion (where $N$ denotes the number of unknowns).

However, a multigrid approach leads to a dilemma when considered in a parallel execution context. While efficiency dictates the coarsest base mesh possible, i.e. in order to accommodate a direct solver, the subsequent small problem size impedes load balancing at large concurrencies~\cite{bastian1998load} and thence leads to a strong degradation of the global parallel efficiency.

Furthermore, most of the available multigrid methods in the scientific literature are limited to second order discretization \cite{brown_multigrid_2005, teunissen_geometric_2019, tomida_athena_2023} with the notable exception of \cite{deriaz_high-order_2023} which have implemented a high-order multigrid solver using shared memory parallelism. While high-order methods are more expensive to use, they allow for a sharper treatment of the solution, thus reducing the amount of needed unknowns to reach a given accuracy (and, in consequence, its computational cost). Such solvers are often sought after when all scales of the problem are being numerically resolved, as evidenced by immersed interfaces methods \cite{gabbard_high-order_2024}.

This work therefore proposes the following combination for the efficient solution of unbounded Poisson problems on adaptive meshes: 
\begin{itemize}
    \item a multigrid solver, which can achieve high-order accuracy,
    \item a multiresolution block-based adaptive grid, which can maintain high-order accuracy at  resolution jumps, 
    \item a FFT-based solver, which solves the coarsest level directly, at a problem size large enough to maintain load balancing.
\end{itemize}

We discuss the resulting implementation, starting with the three-dimensional multiresolution framework, in \sref{section:framework}. 
Then, in \sref{section:methodology}, we show how to use a high-order compact stencil discretization of the Laplacian operator to minimize the performance overheads typically tied to high order methods. We discuss the coupling of the multigrid and FFT solvers, with a particular attention to maintaining load balancing (and parallel performance) and accommodating a variety of boundary conditions including unbounded, symmetric-unbounded and periodic-unbounded.
Validation and performance results are presented in \sref{section:verification} and \sref{section:performances}, demonstrating the accuracy and scalability of the solver. We then close this article with our conclusions and perspectives in \sref{section:conclusion}.

%!TEX root = umg_1_intro.tex
%!TEX encoding = UTF-8 Unicode

\section{Framework}\label{section:framework}
This section provides short introductions to the two existing computational frameworks, which this work builds upon. 
%implemented in \code{C++}. 
%This section provides short introductions A short description of their capabilities is given hereunder.

\subsection{Adaptive multiresolution spatial discretization}
\label{sec:amr_disc}
The adaptive multigrid solver presented hereunder is implemented within \murphy~\cite{gillis_murphy:2022, murphy}, a three-dimensional adaptive multiresolution framework for the solution of PDEs on exascale-era supercomputers. It relies on a collocated block-based domain decomposition stored as a forest of octrees managed using \pforest~\cite{p4est}. 

\murphy\ is constructed from the following tenet: each grid block is a fully independent entity on which finite difference stencils and operators can be applied. In the context of this work, each block is set to a size of $24^3$.
%The core philosophy of \murphy is that each grid block is a fully independent entity on which finite difference stencils and operators can be applied.
Field information updates from adjacent blocks are handled through data mapping, or \textit{ghosting}, following the halo-exchange pattern for distributed memory parallelism. Interpolation at resolution jumps can either be handled using polynomial interpolation or interpolating wavelets of appropriate order. In the present work, this is set to $M+2$ where $M$ is the order of the discretization used. Furthermore, unless explicitly mentioned otherwise, polynomial interpolation is used. 

Grid adaptation is based on the wavelet transform and associated multiresolution analysis. They provide an accurate estimate of the local error committed, which allows for direct and automated error control via a user-specified parameter. Compression of the field information (and the associated error) can then be controlled directly via $\epsilon_r$ and $\epsilon_c$, the refinement and coarsening tolerances, respectively. A grid block is refined into $8$ finer blocks whenever its wavelet transform yields a so-called detail coefficient $\gamma$ above $\epsilon_r$. Alternatively, when all $8$ blocks involve detail coefficients below the coarsening threshold $\epsilon_c$, they revert back to a single grid block. A last point of attention is the so-called 2:1 constraint : two adjacent blocks cannot have a resolution jump larger than one level. Further details on the \murphy\ framework and its implementation can be found at \cite{gillis_murphy:2022}.

\subsection{Fast Fourier Transform-based Poisson solver} There are plenty of FFT-based Poisson solvers such as \texttt{AccFFT} \cite{gholami_accfft_2016}, \texttt{HeFFTE} \cite{krzhizhanovskaya_heffte_2020}, \texttt{P3DFFT} \cite{pekurovsky_p3dfft_2012}, or the recently released open-source library \flups~\cite{caprace_flups_2021, balty_flups_2023, flups}.

This work relies on \flups\ for the coarse-grid direct solver. It indeed features performance optimizations for distributed-memory architectures as well as flexibility with respect to the boundary conditions to be enforced.
The former allows to maintain the scalability of the solver as a whole up to large partitions, while the latter enables the handling of several boundary conditions: homogeneous Dirichlet or Neumann, periodic and unbounded.

%However, \flups is the only solver to have been specifically optimized both for efficient performance on distributed-memory architectures and boundary condition flexibility \cite{caprace_flups_2021, balty_flups_2023}. Thus, we rely on the open-source FFT-based Poisson solver library \flups for this component.

%For this component solver, we rely on \flups, an open-source library that uses FFT-based methods to solve the Poisson equation for any given combination of boundary conditions (unbounded, periodic, homogeneous Dirichlet and Neumann). While there exists alternative FFT-based solver \cite{gholami_accfft_2016, krzhizhanovskaya_heffte_2020, dalcin_fast_2019, pekurovsky_p3dfft_2012}, \flups was specially optimized to perform efficiently on distributed-memory architectures \cite{balty_flups_2023} and boundary condition flexibility \cite{caprace_flups_2021}.

%In this work, we use \flups as the coarse-grid direct solver. As described in the following sections, this is done so to improve the scaling of the multigrid solver as a whole, as well as enable the use of any boundary conditions, including unbounded domains.

Let us recall that FFT-based methods offer best-in-class performances when one considers smooth source terms $f$ on uniform rectilinear grids. 
For problems with unbounded boundary conditions, they essentially perform a convolution of the source term with Green's function, $G(\mathbf{x})$, i.e the solution of the Poisson equation with a Dirac delta function as source term 
\begin{equation}
    \nabla^2G(\mathbf{x}) = \delta(\mathbf{x}),
\end{equation}
and appropriate boundary conditions. The solution for a generic source term then reads
\begin{equation}
    u(\mathbf{x}) = ( G \ast f )(\mathbf{x}),
\end{equation}
or, in spectral space,
\begin{equation}
    \hat{u}(\mathbf{k}) = \hat{G}(\mathbf{k})\, \hat{f}(\mathbf{k}),
\end{equation} 
the starting point of our FFT-based approach.

In practice, we perform the above convolution with the lattice Green's function $G[\mathbf{n}]$ associated to a finite difference operator, i.e. the solution to
\begin{equation}
    \mathcal{L}G[\mathbf{n}] = \delta[\mathbf{n}].
\end{equation}
where $\mathcal{L}$ is a finite difference operator.

The main challenge of this method lies in the computation of the adequate lattice Green function for any combination of boundary conditions \cite{berger1958use, gabbard_lattice_2023, gillman_fast_2010}. Periodic directions are straightforward as they only require the spectral expression of the lattice Green function, which is easily obtained for difference operators. Odd and even boundary conditions, representing homogeneous Dirichlet and Neumann respectively, are just as straightforward but differ from periodic in that they require the use of the appropriate sine or cosine transform instead of the standard discrete Fourier transform \cite{caprace_flups_2021}.

When considering infinite domains where we assume $u \longrightarrow 0$ when $x \longrightarrow \infty$, it is possible to obtain the exact solution to the Poisson equation for compact source terms $f$ using FFT methods. In this case, the numerical domain is a finite region that encloses the sources; one then enforces so-called unbounded boundary conditions. The computation of the associated lattice Green function is not straightforward as it involves the evaluation of singular integrals. As of this writing, such functions have been obtained for three-dimensional domains with one or three unbounded directions \cite{gabbard_lattice_2023}; the integrals involved in the case of two unbounded directions have thus far proved too challenging.

Finally, we also mention the possibility of semi-unbounded configurations. In such cases, an odd or even condition is enforced in a first direction, e.g. $\hat{\boldsymbol{e}}_x$ while an unbounded one is imposed in the opposite, $-\hat{\boldsymbol{e}}_x$. The implementation of such configurations can be achieved by using adequate symmetries, as detailed in \cite{caprace_flups_2021}.
\section{Methodology}\label{section:methodology}

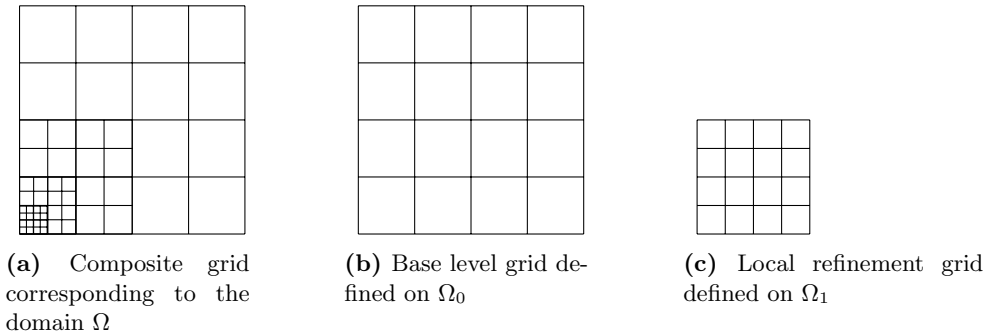
\begin{figure}[h!]
    \centering
    \subfloat[Composite grid corresponding to the domain $\Omega$]{
       \begin{tikzpicture}[scale=1.5]
           \draw (0, 0) grid[step=0.5] (2,2);
           \draw (0, 0) grid[step=0.25] (1,1);
           \draw (0, 0) grid[step=0.125] (0.5,0.5);
           \draw (0, 0) grid[step=0.06125] (0.25,0.25);
       \end{tikzpicture} 
    }
    \hfill
    \subfloat[Base level grid defined on $\Omega_0$]{
       \begin{tikzpicture}[scale=1.5]
            \draw (0, 0) grid[step=0.5] (2,2);
        \end{tikzpicture} 
    }
    \hfill
    \subfloat[Local refinement grid defined on $\Omega_1$]{
       \begin{tikzpicture}[scale=1.5]
           % \draw [very thin, gray] (0, 0) grid[step=0.25] (1.75,1.75);
           \draw (0, 0) grid[step=0.25] (1,1);
           \draw (2.5, 2);
       \end{tikzpicture} 
    }
    \hfill
    \caption{Global AMR grid and its corresponding uniform grid hierarchy}
    \label{fig:grid_hierarchy}
\end{figure}

The core principle of multigrid methods, as discussed in \sref{section:intro}, lies in the accelerated smoothing of the high wavenumber modes of the error, and as such, it is not incompatible with an adaptive discretization, as explained by \cite{trottenberg_multigrid_2000}. Such a framework indeed allows to delimit error contributions both in wavenumber and space. The stack of increasingly coarse grids of classical multigrid methods then gives way to an adaptive multiresolution grid.
The latter is a hierarchy of locally refined grids $\Omega_k$ which, put together, define a composite grid on $\Omega$ as shown at \fref{fig:grid_hierarchy}.
In the context of the present multigrid solver, $\Omega_0$ is assumed to be the finest uniform grid which covers the full domain $\Omega$.

\subsection{Adaptive multigrid cycle}
\begin{algorithm}
    \SetKwComment{Comment}{\,}{\,}
    \newcommand\mygets{\overset{\,k}{\gets}}
    \newcommand{\nosemic}{\renewcommand{\@endalgocfline}{\relax}}
    \newcommand{\dosemic}{\renewcommand{\@endalgocfline}{\algocf@endline}}

    The cycle takes an initial solution $u$ and a source term $f$ and returns an updated solution in-place. Fields are defined on the composite domain $\Omega$. Their local definition on the local domain $\Omega_k$ are noted with the subscript $k$, ranging from $0$ to $\ell$. 
    
    The special assignment $\mygets$ denotes an operation done on the subdomain of the assignee corresponding to $\Omega_k$.\\

    \nosemic
    \For(\Comment*[f]{From fine to coarse}){$k = \ell$ ... $1$}{ 
        \For{$i = 1$ ... $\eta_1$}{
            $u_k \gets \text{smooth}(u_k, f_k)$ \Comment*[r]{Apply pre-smoothing steps}
        }
        $r_k \gets f_k - \nabla^2 u_k$ \Comment*[r]{Compute the residuals}
        $u_{k-1} \mygets \hat{I}^{k-1}_{k}u_k$ \Comment*[r]{Restrict the solution (injection)}
        $r_{k-1} \mygets I^{k-1}_{k}r_k$ \Comment*[r]{Restrict the residuals (full weighting)}
        $f_{k - 1} \mygets r_{k-1} + \nabla^2u_{k-1}$ \Comment*[r]{Compute the right-hand side on $\Omega_k$}
        $\hat{u}_{k-1} \mygets u_{k-1}$ \Comment*[r]{Store the current approximation}
    }
    $u_0 \gets \text{solve}(u_0, f_0)$ \Comment*[r]{Solve on the coarse uniform grid}
    \For(\Comment*[f]{From coarse to fine}){$k = 1$ ... $\ell$}{
        $\epsilon_{k-1} \mygets \hat{u}_{k-1} - u_{k-1}$ \Comment*[r]{Get the correction}
        $u_k \gets u_k + I^{k}_{k-1}\epsilon_{k-1}$ \Comment*[r]{Prolongate the correction}
        \For{$i = 1$ ... $\eta_2$}{
            $u_k \gets \text{smooth}(u_k, f_k)$ \Comment*[r]{Apply post-smoothing steps}
        }
    }
    \caption{Adaptive multigrid V($\eta_1$, $\eta_2$)-cycle}
    \label{alg:amg_vcycle}
\end{algorithm}

The resolution of the Poisson equation consists in the application of successive multigrid cycles to the hierarchy of grids until convergence is reached (typically based on the infinite-norm of the residual field). The presented solver implements the standard V-cycle where one goes down the hierarchy up to the coarsest grid before transferring the correction back up to the finest grids, following the pattern depicted on \fref{fig:amg_vcycle}.

\begin{figure}[h!]
    \centering
    \begin{tikzpicture}
    % Composite grid
    \draw [] (5, 5) node[] {Composite domain $\Omega$};
    \foreach \r in {1,...,4}{
        \begin{scope}[yshift=4.5cm]
            \begin{scope}[yscale = 0.4,xslant=0.4]
                \draw [] (0,0) grid[step=1/(2^\r)] (4 / 2^\r, 4 / 2^\r);
            \end{scope}
        \end{scope}
    }
    % Local grids
    \foreach \r in {1,...,4}{
        \begin{scope}[yshift=\r*1cm - 0.25cm] % y-spacing of the grid
            \begin{scope}[yscale = 0.4,xslant=0.4] % adding the perspective
                \draw (0,0) grid[step=1/(2^\r)] (4 / 2^\r, 4 / 2^\r); % grid with the definition
            \end{scope}
        \end{scope}
    }
    \begin{scope}[xshift=9cm, yshift=1cm]
        % V-cycle
        \draw (0,0) node[]{$\circ$};
        \foreach \r [evaluate=\r as \reval using int(\r - 1)] in {1,...,3} {
            \draw [
              -stealth,
              shorten <=.05cm + \pgflinewidth,
              shorten >=.04cm + \pgflinewidth,
            ]( 0.5*\r - 0.5, \r - 1) -- ( 0.5*\r,\r) node[] {$\bullet$} node[midway, right]{$I^\r_{\reval}$};
            \draw [
              stealth-,
              shorten <=.06cm + \pgflinewidth,
              shorten >=.05cm + \pgflinewidth,
            ](-0.5*\r + 0.5, \r - 1) -- (-0.5*\r,\r) node[] {$\bullet$} node[midway, left ]{$I^{\reval}_\r$};
        }
        % Level labels
        \foreach \r in {0,...,3} {
            %\draw (-4, \r) node[] {level $\r$};
            \draw (-4, \r) node[] {Level $\Omega_\r$};
        }
    \end{scope}
\end{tikzpicture}
    \caption{Adaptive multigrid V-cycle applied to the grid depicted on \fref{fig:grid_hierarchy}}
    \label{fig:amg_vcycle}
\end{figure}

Traditionally, multigrid solvers only transfer the residual field down the hierarchy. However, as the base level $\Omega_0$ is not wholly covered by the local refinement $\Omega_1$, we need to work with the full solution rather than just a correction. We thus have recourse to the so-called full-approximation scheme where both the residual and the solution field are transferred down the hierarchy~\cite{trottenberg_multigrid_2000}. This also helps circumvent any non-linearity issues that may arise at resolution jumps because of the ghost reconstruction schemes which use non-linear polynomial interpolation. 
%Currently, the solver supports the standard V-cycle (as depicted for a 2D grid in \figurename~\ref{fig:amg_vcycle}).

Let us complete the definition of the multigrid solution procedure and introduce some classical operators,
\begin{itemize}
    \item[-] $r = f - \nabla^2 u$, the residual (of the Poisson equation \eqqref{eq:poisson_problem}),
    \item[-] $g_{k-1} = I_{k}^{k-1} g_{k}$, the restriction operator,
    \item[-] $g_{k} = I_{k-1}^{k} g_{k-1}$, the prolongation operator,
    \item[-] $u_{k} = \textrm{smooth}(u_k,f_k)$, the smoothing operator.
\end{itemize}
The smoothing operator is an iteration of the Gauss-Seidel method.
The restriction and prolongation are the inter-grid transfer operators for a field $g_k$ which can be either the residual $r_k$ or the solution $u_k$. For both fields, the prolongation is implemented as a block-wise trilinear interpolation. 
The restriction operator however is specific to the considered field. 
The restriction of the residual $r$ is handled by the three-dimensional full weighting operator; it can be constructed as the tensor product of the one-dimensional operator
\begin{equation}
    I^{k-1}_k = \frac{1}{4}\begin{bmatrix}
        1 & 2 & 1
    \end{bmatrix}.
\end{equation}
The restriction of the solution field $u$ uses the injection operator instead: the coarse value is set to that of the corresponding fine-grid value without any averaging or interpolation. This ensures consistency between the coarse-grid and fine-grid equations everywhere including at the composite grid resolution jumps. The injection operator is simply defined as:
\begin{equation}
    \hat{I}^{k-1}_{k} = \begin{bmatrix}
        1
    \end{bmatrix}
\end{equation}

Considering a hierarchy of locally refined grids such as \fref{fig:amg_vcycle} with levels $k = 0, 1, ..., \ell$, the adaptive multigrid V-cycle is given by \algref{alg:amg_vcycle}.

Once the uniform grid at level $k = 0$ is reached, the solution is computed using a direct solver. The targeted massively parallel applications and the efficiency of the FFT-based direct solver motivate to maintain a discretization fine enough at this level $0$. This coarse problem size is therefore as large as possible, which contributes to preserve load balancing and computational efficiency.
%This is more efficient than coarsening to further levels both from a pure time-to-solution and load balancing perspective.

For the reader interested in more details about adaptive multigrid and multigrid methods in general, we refer to~\cite{trottenberg_multigrid_2000}.

\subsection{Poisson equation discretization}

\begin{figure}[h!]
    \centering
    \begin{tikzpicture}
    % Style stuff
    \tikzstyle{every node}=[font=\footnotesize]
    
    % 7-point Laplacian stencil
    \begin{scope}[yshift=7cm]
        % Scaling factor~
        \draw (-5.5, 0.2, 0) node[]{\normalsize $\frac{1}{h^2}\times$};
        % Connecting line
        \draw[dashed] (-2.5,0,1.5) node[]{$\bullet$} node[below left]{$1$} -- (2.5,0,-1.5) node[]{$\bullet$} node[above right]{$1$};
        % Center point
        \draw (0, 0, 0) node[]{$\bullet$} node[above left]{$-6$};
        % Center cross
        \draw (-1, 0, 0) node[]{$\bullet$} node[above left]{$1$} -- (1, 0, 0) node[]{$\bullet$} node[below right]{$1$};
        \draw (0, -1, 0) node[]{$\bullet$} node[above left]{$1$} -- (0, 1, 0) node[]{$\bullet$} node[above left]{$1$};
    \end{scope}

    % 19-point Laplacian stencil
    \begin{scope}[yshift=3.5cm]
        % Scaling factor~
        \draw (-5.5, 0.2, 0) node[]{\normalsize $\frac{1}{6h^2}\times$};
        % Connecting line
        \draw[dashed] (-2.5,0,1.5) node[]{$\bullet$} node[above left]{$2$} -- (2.5,0,-1.5) node[]{$\bullet$} node[above left]{$2$};
        % Center point
        \draw (0, 0, 0) node[]{$\bullet$} node[above left]{$-24$};
        % Center cross
        \draw (-1, 0, 0) node[]{$\bullet$} node[above left]{$2$} -- (1, 0, 0) node[]{$\bullet$} node[above left]{$2$};
        \draw (0, -1, 0) node[]{$\bullet$} node[above left]{$2$} -- (0, 1, 0) node[]{$\bullet$} node[above left]{$2$};
        % Center square
        \draw (-1, -1, 0) node[]{$\bullet$} node[above left]{$1$} -- (-1, 1, 0) node[]{$\bullet$} node[above left]{$1$} -- (1, 1, 0) node[]{$\bullet$} node[above left]{$1$} -- (1, -1, 0) node[]{$\bullet$} node[above left]{$1$} -- (-1, -1, 0);
        % Front cross
        \draw (-3.5, 0, 1.5) node[]{$\bullet$} node[above left]{$1$} -- (-1.5, 0, 1.5) node[]{$\bullet$} node[above left]{$1$};
        \draw (-2.5, 1, 1.5) node[]{$\bullet$} node[above left]{$1$} -- (-2.5, -1, 1.5) node[]{$\bullet$} node[above left]{$1$};
        % Back cross
        \draw (3.5, 0, -1.5) node[]{$\bullet$} node[above left]{$1$} -- (1.5, 0, -1.5) node[]{$\bullet$} node[above left]{$1$};
        \draw (2.5, 1, -1.5) node[]{$\bullet$} node[above left]{$1$} -- (2.5, -1, -1.5) node[]{$\bullet$} node[above left]{$1$};
    \end{scope}

    % 27-point Laplacian stencil
    \begin{scope}[yshift=0cm]
        % Scaling factor~
        \draw (-5.5, 0.2, 0) node[]{\normalsize $\frac{1}{30h^2}\times$};
        % Connecting line
        \draw[dashed] (-2.5,0,1.5) node[]{$\bullet$} node[above left]{$14$} -- (2.5,0,-1.5) node[]{$\bullet$} node[above left]{$14$};
        % Center point
        \draw (0, 0, 0) node[]{$\bullet$} node[above left]{$-128$};
        % Center cross
        \draw (-1, 0, 0) node[]{$\bullet$} node[above left]{$14$} -- (1, 0, 0) node[]{$\bullet$} node[above left]{$14$};
        \draw (0, -1, 0) node[]{$\bullet$} node[above left]{$14$} -- (0, 1, 0) node[]{$\bullet$} node[above left]{$14$};
        % Center square
        \draw (-1, -1, 0) node[]{$\bullet$} node[above left]{$3$} -- (-1, 1, 0) node[]{$\bullet$} node[above left]{$3$} -- (1, 1, 0) node[]{$\bullet$} node[above left]{$3$} -- (1, -1, 0) node[]{$\bullet$} node[above left]{$3$} -- (-1, -1, 0);
        % Front cross
        \draw (-3.5, 0, 1.5) node[]{$\bullet$} node[above left]{$3$} -- (-1.5, 0, 1.5) node[]{$\bullet$} node[above left]{$3$};
        \draw (-2.5, 1, 1.5) node[]{$\bullet$} node[above left]{$3$} -- (-2.5, -1, 1.5) node[]{$\bullet$} node[above left]{$3$};
        % Front square
        \draw (-3.5, -1, 1.5) node[]{$\bullet$} node[above left]{$1$} -- (-3.5, 1, 1.5) node[]{$\bullet$} node[above left]{$1$} -- (-1.5, 1, 1.5) node[]{$\bullet$} node[above left]{$1$} -- (-1.5, -1, 1.5) node[]{$\bullet$} node[above left]{$1$} -- (-3.5, -1, 1.5);
        % Back cross
        \draw (3.5, 0, -1.5) node[]{$\bullet$} node[above left]{$3$} -- (1.5, 0, -1.5) node[]{$\bullet$} node[above left]{$3$};
        \draw (2.5, 1, -1.5) node[]{$\bullet$} node[above left]{$3$} -- (2.5, -1, -1.5) node[]{$\bullet$} node[above left]{$3$};
        % Back square
        \draw (1.5, -1, -1.5) node[]{$\bullet$} node[above left]{$1$} -- (1.5, 1, -1.5) node[]{$\bullet$} node[above left]{$1$} -- (3.5, 1, -1.5) node[]{$\bullet$} node[above left]{$1$} -- (3.5, -1, -1.5) node[]{$\bullet$} node[above left]{$1$} -- (1.5, -1, -1.5);
    \end{scope}
    
    % % Axis
    \begin{scope}[shift={(-6.5,-1.5)}, scale=0.4]
        \draw[](0,0,0) -- (1,0,0)node[right]{$x$};
        \draw[](0,0,0) -- (0,1,0)node[above]{$y$};
        \draw[dashed](0,0,0) -- (0,0.5,2)node[below left]{$z$};
    \end{scope}
\end{tikzpicture}
    \caption{Left-hand stencils $\Delta^M_h$ for order 2, 4 and 6 respectively (from top to bottom).}
    \label{fig:stencils}
\end{figure}
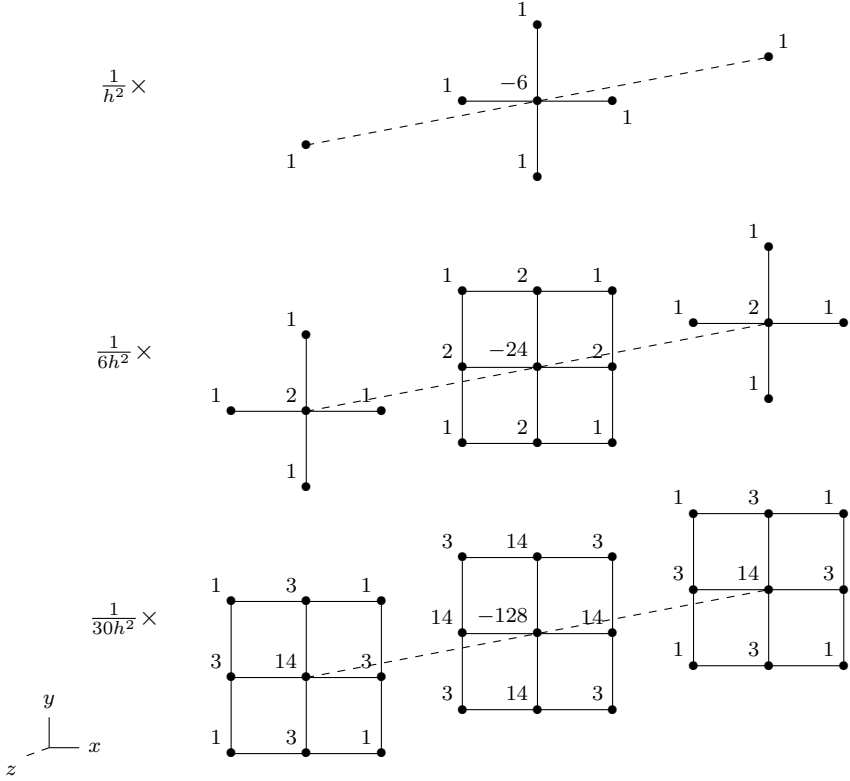

The finite difference discretization of the Laplace operator can be done using cross-shaped stencils, which sum purely one-dimensional stencils. For second order, such a cross-shaped stencil is then given by
\begin{equation}
    \frac{1}{h^2}\left[\delta^2_x + \delta^2_y + \delta^2_z\right]u = f + \mathcal{O}(h^2), \label{eq:cross_discr2}
\end{equation}
where $\delta^2_{x,y,z}$ is the second order central difference operator along the considered direction. It corresponds to
\begin{equation}
    \delta_{x,y,z}^2 =
    \begin{bmatrix}
        1 & -2 & 1
    \end{bmatrix}.
\end{equation}
A three-dimensional representation of the second-order cross Laplacian stencil \eqqref{eq:cross_discr2} is shown in \fref{fig:stencils}. This is the discretization used for the second order-acccurate adaptive multigrid solver.

Higher order discretizations can be achieved using the cross-shaped scheme but they require increasing the width of the stencil. This in turn increases the needed ghost points and associated ghosting communications. 
Such negative impacts on performance can be mitigated if we instead rely on so-called high-order compact \textit{Mehrstellenverfahren} stencils for fourth and sixth order discretizations \cite{spotz_high-order_1996}. These discretizations have a stencil radius of one, thus requiring a single ghost point and are written as:
\begin{equation}
    \Delta^M_hu = R_h^Mf.
\end{equation}
Hence, distinct stencils are applied to the solution $u$ and to the right-hand side $f$. Furthermore, this right-hand side modification can be performed once and one then stores a corrected right-hand side, $f^* = R_h^Mf$ and solve for
\begin{equation}
    \Delta^M_hu = f^*.
\end{equation}

The fourth order Mehrstellen discretization is given by the compact 19-point stencil shown in \fref{fig:stencils} and reads as
\begin{equation}
    \frac{1}{h^2}\left[\delta^2_x + \delta^2_y + \delta^2_z + \frac{1}{6}\left(\delta^2_x\delta^2_y + \delta^2_y\delta^2_z + \delta^2_z\delta^2_x\right) \right]u = f + \frac{1}{12}\left(\delta_x^2 + \delta_y^2 + \delta_z^2\right)f + \mathcal{O}(h^4). \label{eq:mehrstellen_discr4}
\end{equation}

The sixth order stencil involves 27 points; it is shown in \fref{fig:stencils} and reads as
\begin{align}
    \frac{1}{h^2}&\left[\delta^2_x + \delta^2_y + \delta^2_z + \frac{1}{6}\left(\delta^2_x\delta^2_y + \delta^2_y\delta^2_z + \delta^2_z\delta^2_x\right) + \frac{1}{30}\delta^2_x\delta^2_y\delta^2_z\right]u \label{eq:mehrstellen_discr6}\\
    = f &+ \frac{1}{12}\left(\delta_x^2 + \delta_y^2 + \delta_z^2\right)f +\frac{1}{90}\left(\delta^2_x\delta^2_y + \delta_y\delta^2_z + \delta^2_z\delta^2_x\right)f \nonumber \\
    &- \frac{1}{240}\left(\delta_x^4 + \delta_y^4 + \delta_z^4\right)f + \mathcal{O}(h^6). \nonumber
\end{align}
Evidently, due to the fourth order difference operators in the right-hand side, the sixth-order stencil $R^M_h$ has a width of two. Details about the derivation of such stencils can be found at \cite{deriaz_compact_2020}.

\subsection{Unbounded boundary conditions}\label{section:unb}

Unbounded boundary conditions are obtained through the use of \flups\ as the coarse-grid direct solver. This is done under the constraint that every block containing an unbounded boundary must be at the coarsest level of the initial multiresolution grid. Then, the solution on such boundaries is trivially prescribed by the unbounded solution on the coarse grid. This constraint can be theoretically removed by interpolating the so-obtained coarse-block solution to the finer ones, though at the cost of accuracy as shown in \cite{teunissen_geometric_2019}.

The unbounded kernels used by the coarse-grid solver must correspond to the lattice Green's function of the discrete difference operator $\mathcal{L}$ used in the multigrid solver. When the kernel used does not correspond to the multigrid solver's difference operator $\mathcal{L}$, the solver no longer converges even with higher-order kernels. 

To illustrate this compatibility issue, we use a 1D code that follows the same principles as \murphy\ and that implements the above-explained multigrid algorithm on periodic grids. We then perform a convergence test using second-order central finite difference discretization for multigrid-related operations (smoothing, residual, etc.). On the coarse uniform grid, a FFT-based solver is used to recover the solution, using kernels corresponding to second-to-eigth order central finite differences.

\begin{figure}[h!]
    \centering
    \subfloat{
        \includegraphics[width=0.48\textwidth]{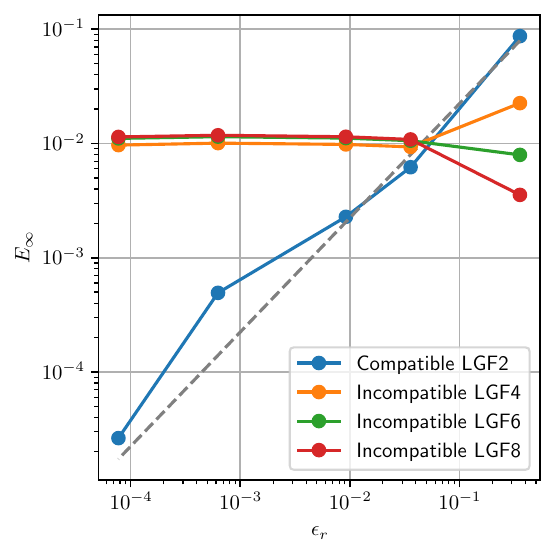}
    }
    \hfill
    \subfloat{
        \includegraphics[width=0.48\textwidth]{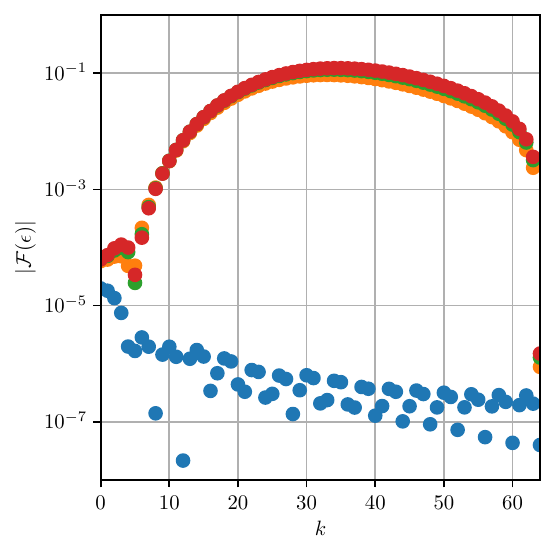}
    }
    \caption{Lattice Green's function (LGF) kernel compatibility issue. Left: convergence of the infinite-norm of the error $E_\infty$ as a function of the refinement criterion $\epsilon_r$. Right: frequency spectrum of the coarse-grid correction.}
    \label{fig:incompatibility}
\end{figure}

The results of this test are shown on \fref{fig:incompatibility}: the multigrids using incompatible kernels (from order 4 to 8) do not converge. This is further highlighted when analyzing the frequency spectrum of the coarse-grid correction: the ones coming from incompatible kernels have a much higher high-frequency content than the one coming from the compatible kernel (second order).

Compatible kernels are known for cases where $1$ or $3$ of the directions contain unbounded boundaries.
\flups\ does provide kernels that work on domains with two unbounded directions, but those kernels do not correspond to the lattice Green's function of the discrete difference operator $\mathcal{L}$ used in the multigrid solver.
However, Spietz et. al \cite{juul_spietz_regularization_2018} suggest to use the following formula,
\begin{equation}
    \widehat{G}_\mathcal{L}(k_x,k_y,k_z) = 
    \begin{cases}
        \widehat{G}_\mathcal{L}^{2D}(k_x, k_y) \hspace{2.18cm}\text{for $k_z = 0$}\\
        1/\sigma_\mathcal{L}(k_x,k_y,k_z) \hspace{1.55cm}\text{elsewhere,}
    \end{cases} \label{eq:2unbkernel}
\end{equation}
where $\widehat{G}_\mathcal{L}^{2D}$ is the already-known LGF for 2D fully unbounded problems \cite{gabbard_lattice_2023} and $\sigma_\mathcal{L}$ is the Fourier symbol of a given discrete difference operator $\mathcal{L}$. and  As shown in the following validation section, this gives us a compatible kernel for domains with unbounded boundaries in two directions.

\subsection{Non-conservation and singular problems}

In \murphy, the interpolation schemes do not guarantee flux continuity across resolution jumps. More explicitly, $\left(\nabla u^h\right)\cdot\hn$ is not enforced to be identical at interfaces between blocks of different resolutions. 

\begin{figure}[h!]
    \centering
   \begin{tikzpicture}[scale=1.5]
       \draw (0, 0) grid[step=0.5] (2,2);
       \draw (0, 0) grid[step=0.25] (1,1);
       \draw (0, 0) grid[step=0.125] (0.5,0.5);
       \draw (0, 0) grid[step=0.06125] (0.25,0.25);

        \draw[red, thick] (0, 1.01) -- (1.01, 1.01) -- (1.01, 0);
        \draw[red] (0, 1.25) node[left]{$\partial\left(\Omega_0 \setminus \Omega_1\right) \setminus \partial\Omega$};
        
        \draw[green!40!gray, thick] (0, 0.99) -- (0.99, 0.99) -- (0.99, 0.99) -- (0.99, 0);
        \draw[green!40!gray, thick] (0, 0.51) -- (0.51, 0.51) -- (0.51, 0.51) -- (0.51, 0);
        \draw[green!40!gray] (0, 0.75) node[left]{$\partial\left(\Omega_1 \setminus \Omega_2\right) \setminus \partial\Omega$};

        \draw[orange, thick] (0, 0.49) -- (0.49, 0.49) -- (0.49, 0.49) -- (0.49, 0);
        \draw[orange, thick] (0, 0.26) -- (0.25, 0.26) -- (0.26, 0.26) -- (0.26, 0);
        \draw[orange] (0, 0.25)  node[left]{$\partial\left(\Omega_2 \setminus \Omega_3\right) \setminus \partial\Omega$};

        \draw[violet, thick] (0, 0.24)  -- (0.24, 0.24) -- (0.24, 0.24) -- (0.24, 0);
        \draw[violet] (0.25, 0) node[below]{$\partial\Omega_3 \setminus \partial\Omega$};

        \draw[blue, thick] (0, 0) grid[step=2] (2, 2) node[above]{$\partial\Omega$};
   \end{tikzpicture} 
    \caption{Resolution jump boundaries $\partial\left(\Omega_{k} \setminus \Omega_{k+1}\right) \setminus \partial\Omega$ of the AMR grid}
    \label{fig:resjump}
\end{figure}
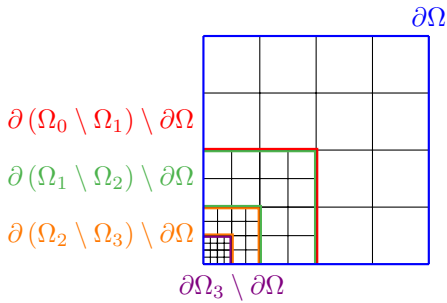

This entails non-conservation on the part of our numerical method.
Indeed, the divergence theorem applied to any region at resolution level $k$, excluding its refined portion, $\Omega_k \setminus \Omega_{k+1}$ (see \fref{fig:resjump}), yields
\begin{equation}
\int_{\Omega_k \setminus \Omega_{k+1}} \nabla^2u^h\,dv = \oint_{\partial(\Omega_k \setminus \Omega_{k+1})} \left(\nabla u^h\right)\cdot\hn \,ds\,.
\end{equation}
where $\hn$ is the outgoing normal.

It is then straightforward to sum over all the levels and identify the contributions from the outer domain boundary $\partial\Omega$ and those from the interfaces between two resolutions
\begin{equation}
    \int_{\Omega}\nabla^2u^h\,dv = \oint_{\partial\Omega} \left(\nabla u^h\right)\cdot\hn \,ds + \sum_{k=0}^\ell \left(\oint_{\partial\left(\Omega_{k} \setminus \Omega_{k+1}\right) \setminus \partial\Omega} \left(\nabla u^h\right)\cdot\hn \,ds \right)\,.
\end{equation}
The second right-hand side term thus constitutes a spurious source term as the normal component of the gradient is discontinuous at resolution interfaces. It can be shown that the resulting surface source term scales as $\mathcal{O}\left(h^{I-1}\right)$ where $I$ is the order of the interpolation scheme used. We recall that $I = M + 2$, where $M$ is the order of the discretization used.

We note that it is generally recommended both for convergence and performance to employ conservative interpolation \cite{trottenberg_multigrid_2000}. However implementing such a scheme in the present collocated (as opposed to cell-centered) multiresolution framework is complex as it depends heavily on the finite difference stencil used. 
The control over the error on the flux discontinuities essentially makes this a non-issue for the considered physical applications and in a local perspective.

There are however consequences for the multigrid solver at a global level. 
Let us consider Poisson problems with only Neumann boundary conditions being imposed, or periodic problems.
It is well-known that the resulting problems are then \textit{singular}: there is a non-empty nullspace, i.e. the solution is non-unique. Furthermore, its existence is contingent on the following conservation property, itself a direct consequence of the divergence theorem,
\begin{equation}
    \int_{\Omega} f^h ds = \oint_{\partial\Omega} g^h \,ds,
    \label{eq:existence_condition}
\end{equation}
where $g^h = \nabla u^h\cdot\hn$ is the prescribed flux on the domain boundary. For our simply connected domain, the nullspace is the constant function and unicity is typically enforced by either arbitrarily setting one nodal value or by requiring the average of the solution to be zero. In our case, we enforce the latter on the coarse-grid solution as suggested by \cite{trottenberg_multigrid_2000}.
The above-discussed non-conservation issue will cause this enforcement to be imperfect and the corresponding solution on the complete grid will have a not-exactly-zero average that depends on the flux discontinuities generated at the resolution jumps. 

Furthermore, in a strict sense, these discontinuities generate an offset in the effective source term which may, in turn, violate the existence condition \eqqref{eq:existence_condition}. This can be highlighted by integrating the residual equation,
\begin{equation}
    \int_\Omega r^h\,ds = \int_{\Omega} f^h ds - \oint_{\partial\Omega} g^h \,ds - \sum_{k=0}^\ell \left(\oint_{\partial\left(\Omega_{k} \setminus \Omega_{k+1}\right) \setminus \partial\Omega} \left(\nabla u^h\right)\cdot\hn \,ds \right).\label{eq:convergence_constraint}
\end{equation}
For a solution to exist, the right-hand side of \eqqref{eq:convergence_constraint} must cancel itself out. 
Because of the discontinuous fluxes, this is not guaranteed in practice. % and leads to the non-convergence of the multigrid solver.

We can also look at this from the point of view of linear algebra: for a discrete linear system $\mathbf{A}\cdot\mathbf{w} = \mathbf{b}$ to have a solution, the right-hand side must be orthogonal to the eigenvector of the transpose of $\mathbf{A}$ corresponding to the null eigenvalue. Denoted $\mathbf{v}$, it must satisfy $\mathbf{v}^T\cdot\mathbf{A} = \mathbf{0}$. Thus, the existence condition \eqqref{eq:existence_condition} can be expressed as
\begin{equation}
    \mathbf{v}^T \cdot \mathbf{b} = \mathbf{0}.
\end{equation}
As previously mentioned, in practice, discretization errors and numerical artifacts, can lead to instances where this condition is not met. The authors of \cite{pozrikidis2001note, trottenberg_multigrid_2000} propose different corrections to apply to the right-hand side $\mathbf{b}$ in order to project it onto the orthogonal complement of the adjoint eigenvector $\mathbf{v}$, thereby rendering the system solvable. However, due to the multiresolution and the complexity of the matrix $\mathbf{A}$ in this case, obtaining the vector $\mathbf{v}$ and the associated correction is non-trivial. Furthermore, as noted by \cite{pozrikidis2001note}, the majority of fluid dynamics application tend to overlook this difficulty. They also mention that the establishment of a consistent numerical scheme that automatically satisfies the existence condition comes at the cost of some loss of generality and greater amount of necessary work.

As the multigrid cycles will not bring the residual $r^n$ to zero due to the above issues, one cannot rely on them to assess convergence of the solution iterates $u^n$.
We then rather use a criterion akin to Cauchy convergence to monitor the iterates,
\begin{equation}
 \frac{\normi{u^{h, n} - u^{h, n+1}}}{\normi{u^{h,n+1}}} < \epsilon.
\end{equation}
This criterion has proven robust and reliable across the battery of tests of the upcoming sections.

\section{Verification}\label{section:verification}
We now consider the maximum local error through 
\begin{equation}
    E_\infty = \left\|u^h - u^{\text{ref}}\right\|_\infty \triangleq \max_{i,j,k,l} \left(\left|u^h_{i,j,k,l}-u^{\text{ref}}_{i,j,k,l}\right|\right)\;,
\end{equation}
where indices run over blocks nodes ($i,j,k$) and grid levels ($l$). This infinity norm conveniently circumvents the need for a quadrature, a tedious operation for the present node-centered discretization. However, the assessment of the convergence behavior for an adaptive grid method still calls for some caution.
The adaptivity indeed precludes the definition of a global grid spacing $h$ and the study of an error behavior $E_\infty(h)$. 

% We recall that 
% \begin{equation}
%     E_\infty = \left\|u^h - u^{\text{ref}}\right\|_\infty \triangleq \max_{i,j,k,\ell} \left(\left|u^h_{i,j,k,\ell}-u^{\text{ref}}_{i,j,k,\ell}\right|\right)\;;
% \end{equation}
% based on the local maximum error, this expression conveniently circumvents the need for a quadrature, a tedious operation for the present node-centered discretization.

We thus study the error behavior against the one parameter available to control the global behavior of the numerical grid, the grid refinement tolerance $\epsilon_r$. Specifically, we adopt a manufactured solution approach with a reference $u^{\text{ref}}(x,y,z)$ and capture this function with an optimally adapted grid. Throughout such a grid, the local wavelet interpolant commits an error which one last refinement step has made smaller than the given $\epsilon_r$. Since the interpolant is a wavelet of order $W$, the local mesh parameter must obey $h^W \propto \epsilon_r$. 

In addition, one should expect our $M$-th order numerical method for the Poisson problem to exhibit a local error $E_\infty$ proportional to $h^M$. The combination of these two asymptotic behaviors then yields the prediction of the maximum local error~\cite{gillis_murphy:2022} with respect to the refinement parameter
\begin{equation}\label{eq:Einf_vs_epsr}
    E_\infty(\epsilon_r) \propto \left(\epsilon_r\right)^{M / W}.
\end{equation}
%
%Studying the convergence behavior of a numerical method on a uniform grid is fairly straightforward: plot the error $E_\infty$ as a function of the grid spacing $h$. However, when considering an adaptive multiresolution grid, we don't have direct control on $h$ nor do we have guarantees that the maximum error $E_\infty$ will be measured at the same physical location when adapting the grid. 
%
%Instead, we choose a reference solution $u_\text{ref}(x,y,z)$ and an upper bound on the error we expect to measure. From these two, an optimal grid given a certain numerical method can be pre-generated. By varying the upper bound on the error and verifying that the solver's numerical error decreases accordingly, we can obtain a proper validation of the numerical method on adaptive multiresolution grids.
%
%In practice, the upper bound on the error depends on the refinement tolerance $\epsilon_r$. Since the grid adaptation is based on wavelets we have that $h \propto \epsilon_r^{1/N}$ where $N$ is the order of the wavelet used. For a given M-th order numerical method, $E_\infty$ is proportional to $h^M$ and we expect the maximum error of such a convergence analysis to behave as \cite{gillis_murphy:2022}
%\begin{equation}
%    E_\infty \propto \left(\epsilon_r\right)^{M / N}.
%\end{equation}
For the remainder of this work, we use the combination $W = M$ unless specified otherwise. This choice arguably provides a balance between the behaviors of the adaptation and the Poisson problem resolution.
%between  a balance  which generates the optimum grid for the given numerical method. 
The following sections present the $E_\infty$ convergence results of the adaptive multigrid solver for different boundary conditions combinations. 

\subsection{Periodic domain}
The reference solution for this case is a Gaussian-like function
\begin{equation}\label{eq:per_uref}
    u_\text{ref}(\mathbf{x}) = \exp\left(-\frac{|\mathbf{x}|^2}{\sigma^2}\right)
\end{equation}
where $\sigma$ is chosen small enough to bring the solution to machine precision-levels on the domain boundary: $u_\text{ref} \ll \epsilon_m$. This essentially ensures that we respect periodic boundary conditions. The right-hand side of the Poisson equation is then
\begin{equation}
    f(\mathbf{x}) = \left(\frac{4|\mathbf{x}|^2}{\sigma^2} - \frac{6}{\sigma^2}\right)\exp\left(-\frac{|\mathbf{x}|^2}{\sigma^2}\right)
\end{equation}

\begin{figure}[h!]
    \centering
    \subfloat{
        \includegraphics[width=0.48\textwidth]{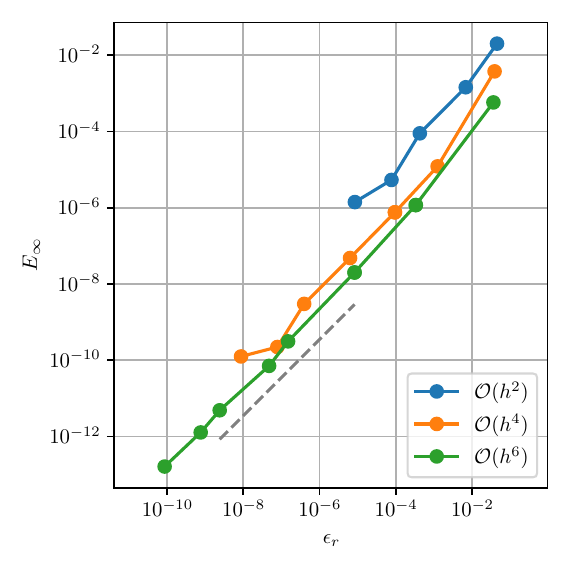}
    }
    \hfill
    \subfloat{
        \includegraphics[width=0.48\textwidth]{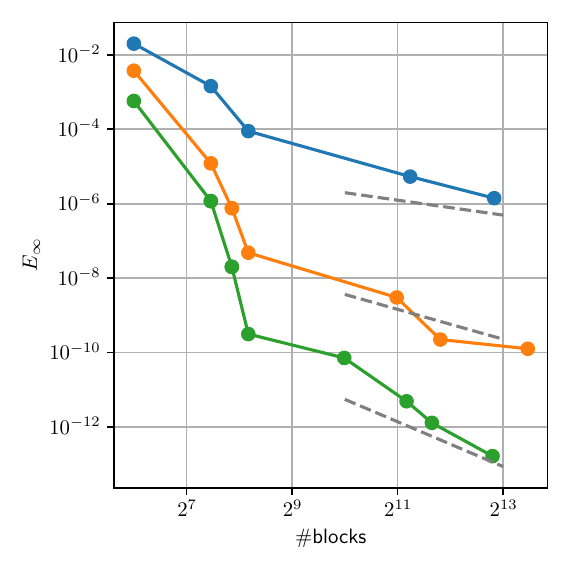}
    }
    \caption{Verification in a periodic unit cube domain for $\sigma = 0.01$ in \eqqref{eq:per_uref}: error behavior $E_\infty$ as a function of the refinement tolerance $\epsilon_r$ (left) and as a function of the resulting number of blocks (right); dashed lines show a linear behavior $E_\infty \propto \epsilon_r$ (left) and the powers $E_\infty \propto (\#blocks)^{-M/3}$ for $M=2,\,4,\,6$ (right, from top to bottom).
%    Multigrid solver convergence for a periodic unit cube domain and $\sigma = 0.01$. The infinite norm of the error $E_\infty$ is shown as a function of the refinement criterion $\epsilon_r$ (left) and the number of blocks (right). Dashed lines represent the expected linear behavior (left) and "equivalent uniform resolution convergence", i.e. $\mathcal{O}(\#blocks^{M/3})$ (right). 
    }
    \label{fig:validation_3per}
\end{figure}

\Fref{fig:validation_3per} sheds some light on the error behavior through the refinement process. Firstly, the error and refinement tolerance can be seen to obey the linear relationship expected from \eqqref{eq:Einf_vs_epsr} and our choice of $M=W$. The error can then be considered from the perspective of the activated degrees-of-freedom, or equivalently the memory footprint of the solver, here through the number of blocks. One can distinguish two phases: at coarse resolutions, the convergence behavior is spectral-like while fine resolutions adopt a polynomial behavior; these two phases are actually expected. The solution is indeed characterized by a compact support and fine features in the middle of the domain: at low resolutions, or equivalently large $\epsilon_r$'s, the adaptation then only adds blocks in the vicinity of this well-identified feature.
More demanding $\epsilon_r$'s lead to a saturation-like effect. This is due to the 2:1 constraint of the multiresolution framework (see \sref{sec:amr_disc}). At some point, this constraint causes the refinement to concern the whole domain and thus leads to a conventional polynomial behavior dictated by the solver order, in a fashion similar to a uniform grid refinement process.

%We also show the error obtained as a function of the number of blocks (each block containing $24^3$ unknowns). From this, we can distinguish two phases: the first part where the solver exhibits spectral-like convergence and the second where we observe a more conventional polynomial convergence (akin to what we would expect from a uniform grid convergence validation). This is because the adaptation process initially adds blocks around a single point only (where the sharp Gaussian-like function is located). Past a certain point, each block added triggers the 2:1 constraint, causing almost all blocks to refine themselves in a similar fashion to a traditional uniform grid convergence study.

% For the sake of completeness, we also show the error obtained as a function of the number of blocks (each block containing $24^3$ unknowns). 
%This showcases the advantages of using a higher-order discretisation: for a given error, it is always cheaper (in term of unknowns) than a lower-ordered one.

\subsection{Unbounded domains} 
For the various combinations of periodic/unbounded boundary conditions, we construct the reference solution as a product of one-dimensional functions,
\begin{equation}
    u_\text{ref}(x, y, z) = u_\text{x}(x) u_\text{y}(y) u_\text{z}(z)
\end{equation}
from which the right-hand side of \eqqref{eq:poisson_problem} is trivially obtained as
\begin{equation}
    f(x, y, z) = \frac{du^2_\text{x}}{dx^2}u_\text{y}u_\text{z} + u_\text{x}\frac{du^2_\text{y}}{dy^2}u_\text{z} + u_\text{x}u_\text{y}\frac{du^2_\text{z}}{dz^2}\;.
\end{equation}
Each component function is designed to agree with the boundary conditions along its direction. We thus have $u_{x,y,z} = u_\text{per}$ or $u_\text{unb}$ for periodic and unbounded directions respectively, with 
\begin{align}
    u_\text{per}(x) &= \exp\left[\sin\left(2\pi x\right)\right] - 1\\
    u_\text{unb}(x) &= \exp\left[10\left(1 - \frac{1}{1 - \left(2x - 1\right)^2}\right)\right],
    % u_\text{per}(x) &= \exp\left[\sin\left(\frac{2\pi f x}{L}\right)\right] - 1\\
    % u_\text{unb}(x) &= \exp\left[10\left(1 - \frac{1}{1 - \left(\frac{2x}{L} - 1\right)^2}\right)\right],
\end{align}
where one notes the regularity and the compact character of the unbounded component.

\begin{figure}[h!]
    \centering
    \subfloat{
        \includegraphics[width=0.48\textwidth]{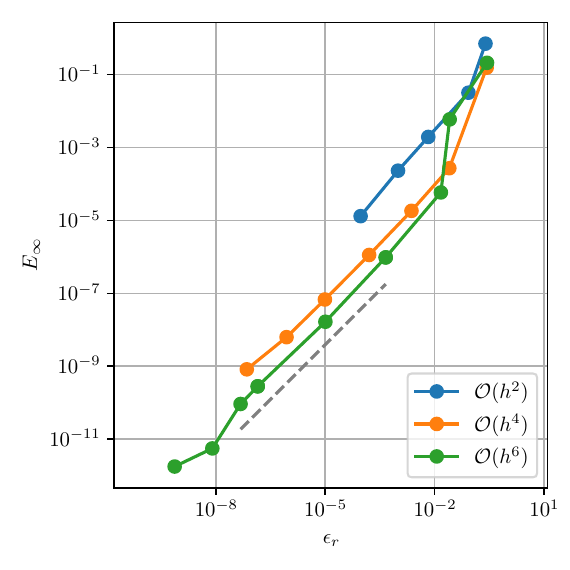}
    }
    \hfill
    \subfloat{
        \includegraphics[width=0.48\textwidth]{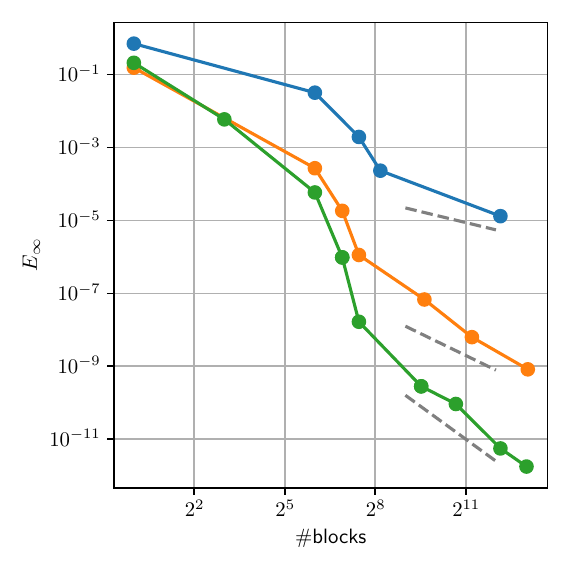}
    }
    \caption{Verification in an unbounded unit cube domain: error behavior $E_\infty$ as a function of the refinement tolerance $\epsilon_r$ (left) and as a function of the resulting number of blocks (right); dashed lines show a linear behavior $E_\infty \propto \epsilon_r$ (left) and the powers $E_\infty \propto (\#blocks)^{M/3}$ for $M=2,\,4,\,6$ (right, from top to bottom).
    % Multigrid solver convergence for a unit cube unbounded domain. The infinite norm of the error $E_\infty$ is shown as a function of the refinement criterion $\epsilon_r$ (left) and the number of blocks (right).
    }
    \label{fig:validation_3unb}
\end{figure}

\begin{figure}[h!]
    \centering
    \subfloat{
        \includegraphics[width=0.48\textwidth]{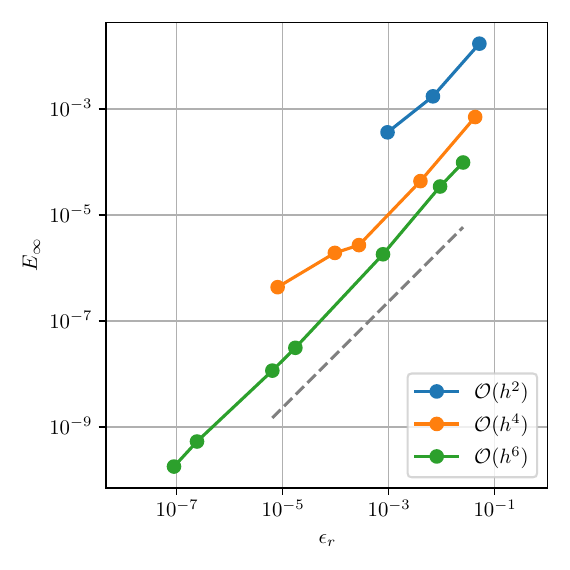}
    }
    \hfill
    \subfloat{
        \includegraphics[width=0.48\textwidth]{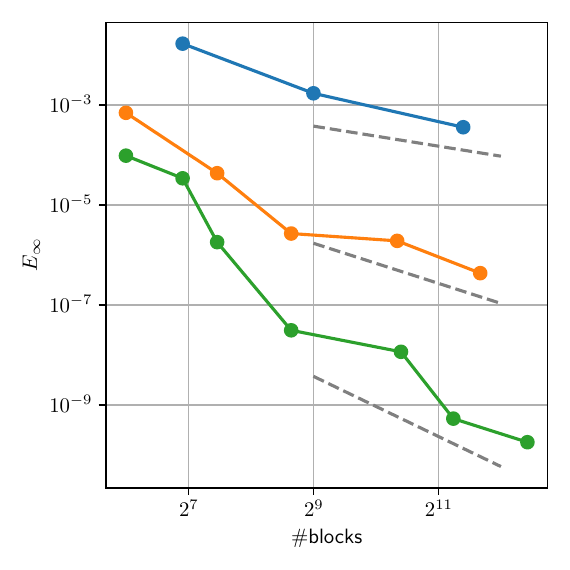}
    }
    \caption{Verification in a 2-unbounded and 1-periodic unit cube domain: error behavior $E_\infty$ as a function of the refinement tolerance $\epsilon_r$ (left) and as a function of the resulting number of blocks (right); dashed lines show a linear behavior $E_\infty \propto \epsilon_r$ (left) and the powers $E_\infty \propto (\#blocks)^{-M/3}$ for $M=2,\,4,\,6$ (right, from top to bottom).
    %Multigrid solver convergence for a unit cube 2-unbounded and 1-periodic domain. The infinite norm of the error $E_\infty$ is shown as a function of the refinement criterion $\epsilon_r$ (left) and the number of blocks (right).
    }
    \label{fig:validation_2unb}
\end{figure}

\begin{figure}[h!]
    \centering
    \subfloat{
        \includegraphics[width=0.48\textwidth]{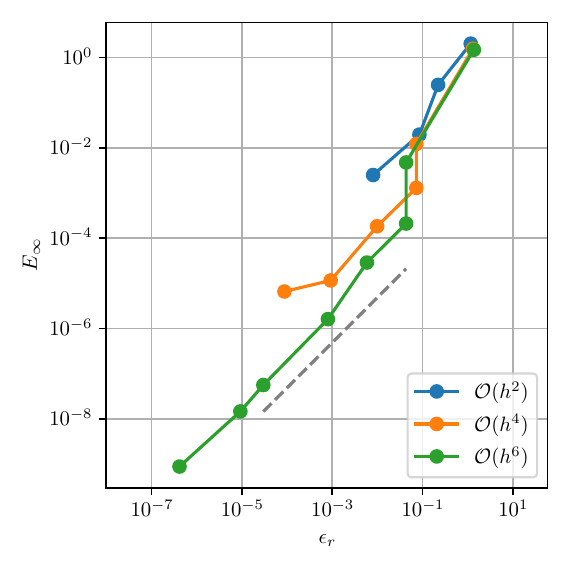}
    }
    \hfill
    \subfloat{
        \includegraphics[width=0.48\textwidth]{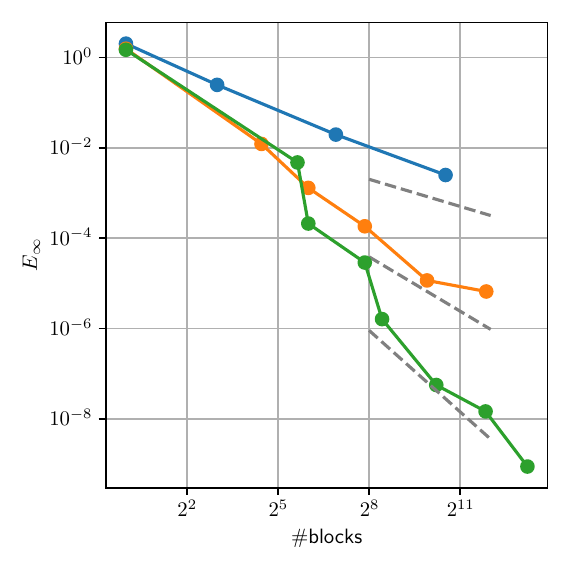}
    }
    \caption{Verification in an 1-unbounded and 2-periodic unit cube domain: error behavior $E_\infty$ as a function of the refinement tolerance $\epsilon_r$ (left) and as a function of the resulting number of blocks (right); dashed lines show a linear behavior $E_\infty \propto \epsilon_r$ (left) and the powers $E_\infty \propto (\#blocks)^{-M/3}$ for $M=2,\,4,\,6$ (right, from top to bottom).
    % Multigrid solver convergence for a unit cube 1-unbounded and 2-periodic domain. The infinite norm of the error $E_\infty$ is shown as a function of the refinement criterion $\epsilon_r$ (left) and the number of blocks (right).
    }
    \label{fig:validation_1unb}
\end{figure}

\Fref{fig:validation_3unb,fig:validation_2unb,fig:validation_1unb} show the convergence results for domains with three, two and one unbounded directions respectively. We recover the expected linear relationship between error and refinement tolerance for the three cases but we can differentiate between their behaviors with respect to the degrees-of-freedom.
For three unbounded directions, the features of the manufactured solution are quite localized and lead to a two-phase behavior similar to that of the fully-periodic test case in \fref{fig:validation_3per}.
Meanwhile, the manufactured solutions for the cases with one and two unbounded directions exhibit features uniformly distributed along the periodic directions. The 2:1 constraint then drives a propagation of adaptation from either a plane or a line for the one and two unbounded directions cases. These cases then transition to the polynomial regime much earlier, making the spectral-like convergence invisible.

Finally, let us make it clear that the two regimes observed for the error against degrees-of-freedom and the transitions in between are quite sensitive to the distribution of features in the source term--and thence solution-- and by the types of boundary conditions. They therefore have to be understood as limiting cases for the error behavior in general.

%already trigger the saturation effect discussed above at coarse refinements as the features are quite distributed in the periodic directions.

%actually compact where the adaptation tends to focus around a single point (the test function in this case is also compact Gaussian-like) and we observe similar phases (spectral-like convergence first followed by polynomial-like) as in the 3-periodic cases. For two and one unbounded directions however, the adaptation tends to develop along a line and a plane (respectively), thus triggering the 2:1 constraint much earlier. Hence, we no longer observe the spectral-like convergence behavior of the 3-periodic and 3-unbounded cases.
%This also further validates the efficiency of the higher-order discretisation compared to its lower-order counterparts.

\section{Parallel performance}\label{section:performances}

The solver performances are tested using the Biot-Savart law, a variation of the standard Poisson problem,
\begin{equation}
    \nabla^2 \boldsymbol{u} = - \nabla\times \boldsymbol{f}
\end{equation}
with vector fields $\boldsymbol{u}$ and $\boldsymbol{f}$.
Such a problem rises in electromagnetism or in incompressible fluid dynamics. In the latter, one often has to recover the flow velocity field $\boldsymbol{u}$ from its curl $\boldsymbol{f}$, typically referred to as vorticity $\boldsymbol{\omega}$.
%In fluid dynamics, it is commonly used to recover the velocity field from a given vorticity field. 
We perform our study on the case of a compact vortex tube aligned in the $z$-direction and centered within unit cube domain. This domain is unbounded in the $x$- and $y$-directions and periodic in the $z$-direction. The expression of the vorticity is,
\begin{equation}
    \boldsymbol{\omega}(x, y, z) = \omega_z(r)\, \boldsymbol{\hat{e}}_z
\end{equation}
where $r = \sqrt{\left(x - 1/2\right)^2 + \left(y-1/2\right)^2}$ and $\omega_z(r)$ is given by
\begin{equation}
    \omega_z(r) = 
    \begin{cases}
        \frac{1}{2\pi}\frac{2}{R^2}\frac{1}{E_2(1)}\exp\left(-\frac{1}{1 - \left(\frac{r}{R}\right)^2}\right) \hspace{1cm}\text{if }r\leq R\\
        0 \hspace{5.2cm} \text{otherwise}
    \end{cases}
\end{equation}
where $R$ is the radius of the vortex tube and $E_2$ is the generalized exponential integral function. The corresponding analytical velocity is known and given by \cite{winckelmans2004encyclopedia}
\begin{equation}
     \boldsymbol{u}(x, y, z) =-\sin\left(\theta\right) u_\theta(r)\, \boldsymbol{\hat{e}}_x + \cos\left(\theta\right)u_\theta(r)\, \boldsymbol{\hat{e}}_y
%    u_\theta(r) = \left\{-\sin\left(\theta\right)u_\theta(r),\,\cos\left(\theta\right)u_\theta(r),\,0\right\}
\end{equation}
where $\theta = \arctan(y,x)$ and $u_\theta(r)$ is given by
\begin{equation}
    u_\theta(r) = 
    \begin{cases}
        \frac{1}{2\pi r} \left[1 - \frac{1}{E_2(1)}\left(1 - \left(\frac{r}{R}\right)^2\right)E_2\left(\frac{1}{1 - \left(\frac{r}{R}\right)^2}\right)\right] \hspace{1cm}\text{if }r\leq R\\
        \frac{1}{2\pi r} \hspace{7.08cm} \text{otherwise}.
    \end{cases}
\end{equation}

In order to be closer to an actual application, the multiresolution grid is here adapted based on the right-hand side, vorticity $\boldsymbol{\omega}$ rather than the analytical solution $\boldsymbol{u}$. A weak scaling study is obtained by successively duplicating the domain along the tube's direction. For this test, the refinement tolerance was set to $\epsilon_r = 10^{-2}$.

This test was performed on the CPU partition of the LUMI supercomputer, a Tier-0 European HPC system. Each LUMI node has two AMD EPYC 7763 CPUs, for a total of $128$ cores per node. The interconnect is a $200$ Gb/s Slingshot-11 network.

\begin{figure}[h!]
    \centering
    \subfloat[Time/cycle]{
        \includegraphics[width=0.48\textwidth]{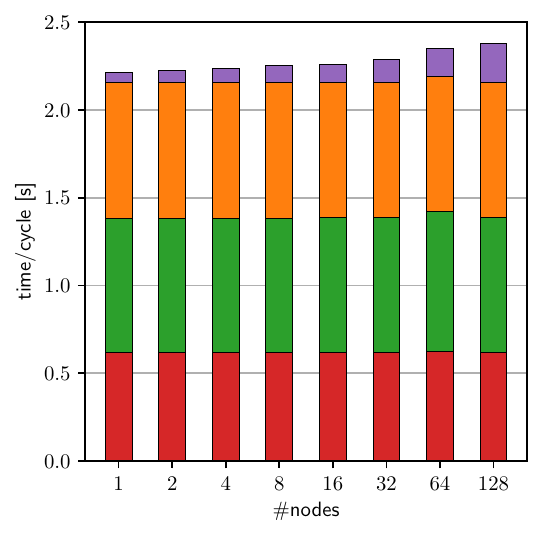}
    }
    \hfill
    \subfloat[Efficiency]{
        \includegraphics[width=0.48\textwidth]{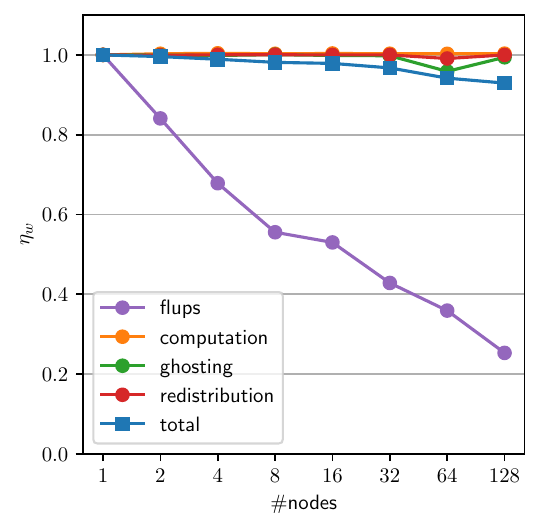}
    }
    \caption{Performance: weak scaling for fourth order stencil ($60$ blocks or $829440$ grid points per rank)}
    \label{fig:weak_scaling_o4}
\end{figure}

\begin{figure}[h!]
    \centering
    \subfloat[Time/cycle]{
        \includegraphics[width=0.48\textwidth]{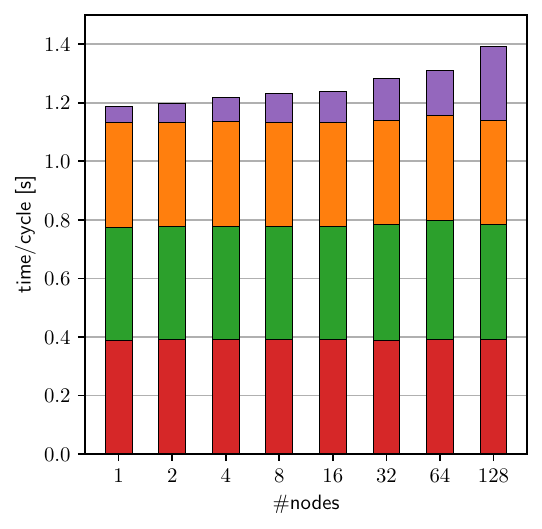}
    }
    \hfill
    \subfloat[Efficiency]{
        \includegraphics[width=0.48\textwidth]{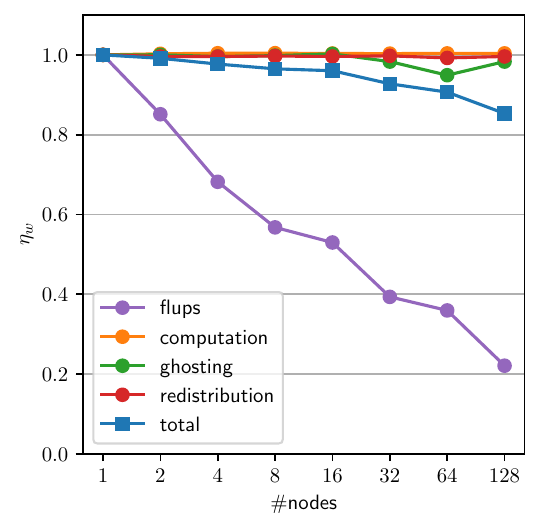}
    }
    \caption{Performance: weak scaling for sixth order stencil ($25$ blocks or $345600$ grid points per rank)}
    \label{fig:weak_scaling_o6}
\end{figure}

Timings for a single V(3,3)-cycle and its weak scaling efficiency $\eta_w$ are shown in \fref{fig:weak_scaling_o4,fig:weak_scaling_o6} for the 4th and 6th order discretization respectively. We divide the computational operations of the V-cycle into four groups: (i) ghosting; (ii) redistribution of blocks among processing units: this has to be performed at each level change to maintain load balancing; (iii) \flups: the FFT-based direct solve on the coarsest grid; and (iv) computations: those are essentially stencil operations, namely smoothing, restriction and prolongation. All these computational tasks exhibit perfect weak scaling efficiency except for \flups. This is however expected due to the low number of unknowns (here, four blocks or $55296$ unknowns) per rank at that level. 

Indeed, over the weak scaling procedure, the problem size in the unbounded directions remains at eight blocks across, i.e. $8\times24=192$ grid points.
At scale for $128$ nodes, the first transform (done along the periodic direction) is partitioned among $128^2$ processing units. Each process has then, on average, $(192/128)^2 = 2.25$ FFT to perform. Obviously, this is very low and actually closer to a stringent strong scaling test.
%over the weak scaling procedure, the problem size in the unbounded directions remains at eight blocks across, i.e. $8\times24=192$. This leads to  Given that there are exactly $192^2$ transforms to execute in that direction, each process has, on average, $1.5$ FFT to do, which is very low.
Ultimately, the multigrid solver in its whole manages to retain an excellent weak scaling efficiency as \flups\ absorbs the brunt of the parallel inefficiency but crucially for a limited contribution to the execution time.

\begin{figure}[h!]
    \centering
    \subfloat[Residual]{\includegraphics[width=0.48\textwidth]{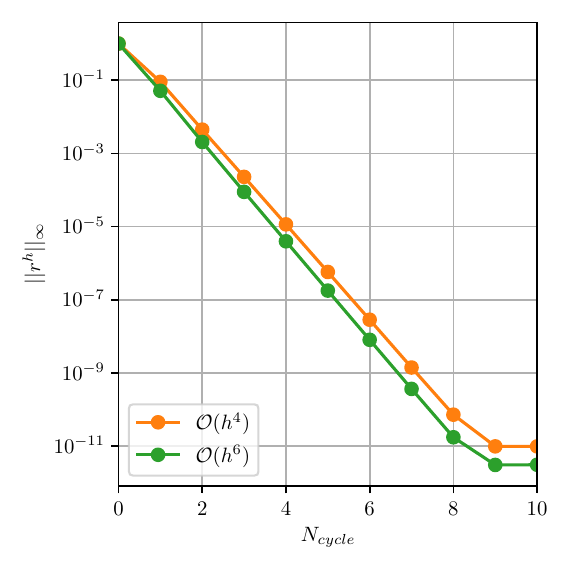}}
    \hfill
    \subfloat[Error]{
        \includegraphics[width=0.48\textwidth]{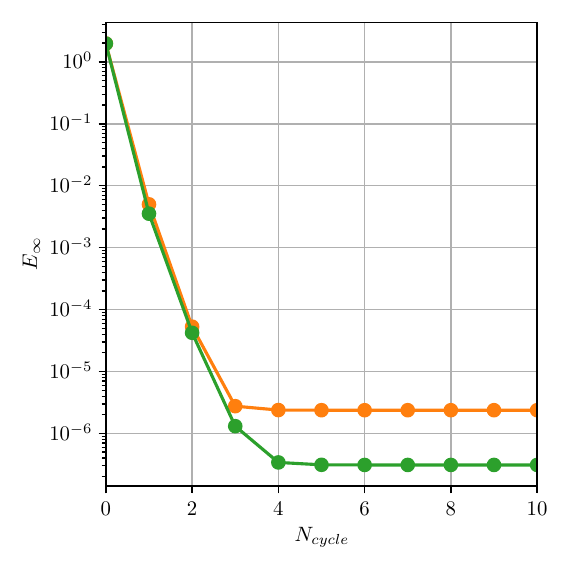}
        \label{subfig:weak_res_convergence}
    }
    \caption{Performance: convergence over cycles}
    \label{fig:weak_convergence}
\end{figure}

For the sake of completeness, \fref{fig:weak_convergence} shows the behaviors of both the residual and the error over the successive cycles. 
One sees that three V-cycles are sufficient to reduce the residual by four orders of magnitude. For the problem at hand, \fref{subfig:weak_res_convergence} indicates that this corresponds to essentially reaching the discretization error. Those choices lead to times-to-solution of approximately $6.5$ and $3.6$ seconds for the fourth and sixth order solvers, respectively. 
In spite of its inherent higher computational cost, the sixth-order discretization here demonstrates its superiority with respect to a time-vs-acccuracy performance.

%!TEX root = umg_6_conclusion.tex
%!TEX encoding = UTF-8 Unicode

\section{Conclusions}\label{section:conclusion}
% Recap
In this paper, we have a introduced and verified a massively parallel Poisson solver for adaptive multiresolution collocated grids. This solver relies on the combination of a multigrid solver with a direct Fourier-based solver for the coarsest level. The former is implemented within the \murphy\ framework while the latter is provided by the \flups\ library, both being open-source efforts.

% Salient points
Beyond its adaptive multiresolution discretization, the presented solution technique stands out by several features and computational aspects:
(i) it can handle periodic, unbounded, and homogeneous Dirichlet or Neumann boundary conditions;
(ii) in case of unbounded domains, the solver enforces consistency between the Finite Difference stencils and Lattice Green's functions, a necessary condition to ensure the accuracy of the complete solver, also at high order; 
(iii) the use of a high-performance direct FFT-based solver on a still large problem at the coarsest level mitigates typical load balancing issues of multigrid solvers;
(iv) the use of compact stencils ensures that no additional communication cost is incurred when using higher-order discretizations.

% results
On simple test-cases, we have found the convergence and accuracy-resolution trade-off to behave as expected for all the possible domain configurations. In particular, the error trend with respect to the discretization size was found to exhibit two regimes. While it does follow the expected spectral-like behavior at coarse resolutions, the $2:1$ refinement constraint accelerates the propagation of the refinement over the whole domain for very fine resolutions, thence leading a uniform resolution-like polynomial behavior. This obviously depends in the distribution of fine scales in the problem but could nevertheless be alleviated by relaxing the $2:1$ constraint.
Parallel performance was found to be excellent for a large problem representative of fluid mechanics applications. This test was most stringent and stressed the load balancing for the coarse level FFT solver, leading to the poor weak scaling of that individual component.
This has little impact on the overall weak scalability as the \flups\ FFT solver has a very limited cost with respect to the remaining components, which do have perfect weak scaling.

% Impact and perspectives
The present solver brings an efficient and flexible solution for a potentially broad public in computational physics. Incidentally, it is already used as a component in the recent implementation of a semi-Lagrangian method for vortex dynamics problems \cite{balty_2025}.
Still, the present solver could be improved by making the refinement more flexible, i.e. the unbounded level and $2:1$ constraints. This is the subject of future work. Finally, let us mention that ongoing developments focus on porting \flups,  \murphy\ and the multigrid solver to heterogeneous CPU-GPU architectures.

\section*{Acknowledgement}
    Computational resources have been provided by the supercomputing facilities of the Université catholique de Louvain (CISM/UCL) and the Consortium des Équipements de Calcul Intensif en Fédération Wallonie Bruxelles (CÉCI) funded by the Fond de la Recherche Scientifique de Belgique (F.R.S.-FNRS) under convention 2.5020.11 and by the Walloon Region. 
    The present research benefited from computational resources made available on Lucia, the Tier-1 supercomputer of the Walloon Region, infrastructure funded by the Walloon Region under the grant agreement n°1910247.
    We acknowledge LUMI-BE for awarding this project access to the LUMI supercomputer, owned by the EuroHPC Joint Undertaking, hosted by CSC (Finland) and the LUMI consortium through a LUMI-BE Regular Access call. LUMI-BE is a joint effort from BELSPO (federal), SPW Économie Emploi Recherche (Wallonia), Department of Economy, Science \& Innovation (Flanders) and Innoviris (Brussels).

\bibliographystyle{siam}
\bibliography{AllPublications}

@article{gillis_murphy:2022,
	title = {{MURPHY}---{A} {Scalable} {Multiresolution} {Framework} for {Scientific} {Computing} on {3D} {Block}-{Structured} {Collocated} {Grids}},
	volume = {44},
	issn = {1064-8275, 1095-7197},
	url = {https://epubs.siam.org/doi/10.1137/21M141676X},
	doi = {10.1137/21M141676X},
	language = {en},
	number = {5},
	urldate = {2023-01-18},
	journal = {SIAM Journal on Scientific Computing},
	author = {Gillis, Thomas and van Rees, Wim M.},
	month = oct,
	year = {2022},
	pages = {C367--C398},
}

@article{p4est,
  author = {Carsten Burstedde and Lucas C. Wilcox and Omar Ghattas},
  title = {{\texttt{p4est}}: Scalable Algorithms for Parallel Adaptive Mesh
           Refinement on Forests of Octrees},
  journal = {SIAM Journal on Scientific Computing},
  volume = {33},
  number = {3},
  pages = {1103-1133},
  year = {2011},
  doi = {10.1137/100791634}
}

@article{spotz_high-order_1996,
	title = {A high-order compact formulation for the {3D} {Poisson} equation},
	volume = {12},
	issn = {1098-2426},
	url = {https://onlinelibrary.wiley.com/doi/abs/10.1002/%28SICI%291098-2426%28199603%2912%3A2%3C235%3A%3AAID-NUM6%3E3.0.CO%3B2-R},
	doi = {10.1002/(SICI)1098-2426(199603)12:2<235::AID-NUM6>3.0.CO;2-R},
	language = {en},
	number = {2},
	urldate = {2023-02-09},
	journal = {Numerical Methods for Partial Differential Equations},
	author = {Spotz, W. F. and Carey, G. F.},
	year = {1996},
	pages = {235--243},
}

@article{deriaz_compact_2020,
	title = {Compact finite difference schemes of arbitrary order for the {Poisson} equation in arbitrary dimensions},
	volume = {60},
	issn = {1572-9125},
	url = {https://doi.org/10.1007/s10543-019-00772-5},
	doi = {10.1007/s10543-019-00772-5},
	language = {en},
	number = {1},
	urldate = {2023-02-17},
	journal = {BIT Numerical Mathematics},
	author = {Deriaz, Erwan},
	month = mar,
	year = {2020},
	pages = {199--233},
}

@book{trottenberg_multigrid_2000,
	title = {Multigrid},
	publisher = {Elsevier},
	author = {Trottenberg, Ulrich and Oosterlee, Cornelius W. and Schuller, Anton},
	year = {2000},
}

@article{teunissen_geometric_2019,
	title = {A geometric multigrid library for quadtree/octree {AMR} grids coupled to {MPI}-{AMRVAC}},
	volume = {245},
	journal = {Computer Physics Communications},
	author = {Teunissen, Jannis and Keppens, Rony},
	year = {2019},
	pages = {106866},
}

@article{tomida_athena_2023,
	title = {The {Athena}++ {Adaptive} {Mesh} {Refinement} {Framework}: {Multigrid} {Solvers} for {Self}-gravity},
	volume = {266},
	issn = {0067-0049, 1538-4365},
	shorttitle = {The {Athena}++ {Adaptive} {Mesh} {Refinement} {Framework}},
	url = {https://iopscience.iop.org/article/10.3847/1538-4365/acc2c0},
	doi = {10.3847/1538-4365/acc2c0},
	language = {en},
	number = {1},
	urldate = {2025-03-05},
	journal = {The Astrophysical Journal Supplement Series},
	author = {Tomida, Kengo and Stone, James M.},
	month = may,
	year = {2023},
	pages = {7},
}

@article{hejlesen_high_2013,
	title = {A high order solver for the unbounded {Poisson} equation},
	volume = {252},
	issn = {0021-9991},
	url = {https://www.sciencedirect.com/science/article/pii/S0021999113004324},
	doi = {10.1016/j.jcp.2013.05.050},
	language = {en},
	urldate = {2022-12-15},
	journal = {Journal of Computational Physics},
	author = {Hejlesen, Mads Mølholm and Rasmussen, Johannes Tophøj and Chatelain, Philippe and Walther, Jens Honoré},
	month = nov,
	year = {2013},
	pages = {458--467},
}

@article{caprace_flups_2021,
	title = {{FLUPS}: {A} {Fourier}-{Based} {Library} of {Unbounded} {Poisson} {Solvers}},
	volume = {43},
	issn = {1064-8275},
	shorttitle = {{FLUPS}},
	url = {https://epubs.siam.org/doi/abs/10.1137/19M1303848},
	doi = {10.1137/19M1303848},
	number = {1},
	urldate = {2022-02-09},
	journal = {SIAM Journal on Scientific Computing},
	author = {Caprace, Denis-Gabriel and Gillis, Thomas and Chatelain, Philippe},
	month = jan,
	year = {2021},
	pages = {C31--C60},
}

@article{balty_flups_2023,
	title = {{FLUPS} - {A} {Flexible} and {Performant} {Massively} {Parallel} {Fourier} {Transform} {Library}},
	volume = {34},
	issn = {1558-2183},
	doi = {10.1109/TPDS.2023.3254302},
	number = {7},
	journal = {IEEE Transactions on Parallel and Distributed Systems},
	author = {Balty, Pierre and Chatelain, Philippe and Gillis, Thomas},
	month = jul,
	year = {2023},
	pages = {2011--2024},
}

@article{chatelain_fourier-based_2010,
	title = {A {Fourier}-based elliptic solver for vortical flows with periodic and unbounded directions},
	volume = {229},
	issn = {0021-9991},
	url = {https://www.sciencedirect.com/science/article/pii/S0021999110000021},
	doi = {10.1016/j.jcp.2009.12.035},
	number = {7},
	urldate = {2024-11-28},
	journal = {Journal of Computational Physics},
	author = {Chatelain, Philippe and Koumoutsakos, Petros},
	month = apr,
	year = {2010},
	pages = {2425--2431},
}

@book{hockney_computer_2021,
	title = {Computer {Simulation} {Using} {Particles}},
	isbn = {978-0-367-80693-4},
	publisher = {Taylor \& Francis},
	author = {Hockney, R. W. and Eastwood, J. W.},
	year = {1988},
	doi = {10.1201/9780367806934},
}

@incollection{chatelain_vortex_2008,
	series = {Lecture {Notes} in {Computer} {Science}},
	title = {Vortex {Methods} for {Massively} {Parallel} {Computer} {Architectures}},
	volume = {5336},
	booktitle = {High {Performance} {Computing} for {Computational} {Science} - {VECPAR} 2008},
	publisher = {Springer Berlin / Heidelberg},
	author = {Chatelain, Philippe and Curioni, Alessandro and Bergdorf, Michael and Rossinelli, Diego and Andreoni, Wanda and Koumoutsakos, Petros},
	editor = {Palma, José and Amestoy, Patrick and Daydé, Michel and Mattoso, Marta and Lopes, João},
	year = {2008},
	pages = {479--489},
}

@misc{ji_fourth_2025,
	title = {A fourth order sharp immersed method for the incompressible {Navier}-{Stokes} equations with stationary and moving boundaries and interfaces},
	url = {http://arxiv.org/abs/2508.15083},
	doi = {10.48550/arXiv.2508.15083},
	language = {en},
	urldate = {2025-10-01},
	publisher = {arXiv},
	author = {Ji, Xinjie and Shen, Changxiao Nigel and Rees, Wim M. van},
	month = aug,
	year = {2025},
	note = {arXiv:2508.15083 [physics]},
}

@misc{winckelmans2004encyclopedia,
  title={Encyclopedia of computational mechanics-Volume 1. Chapter 5: Vortex Methods, chapter 5},
  author={Winckelmans, GS},
  publisher={Wiley Online Library}
}

@misc{gabbard_lattice_2023,
	title = {Lattice {Green}'s {Functions} for {High} {Order} {Finite} {Difference} {Stencils}},
	url = {http://arxiv.org/abs/2309.13503},
	doi = {10.48550/arXiv.2309.13503},
	urldate = {2023-09-26},
	publisher = {arXiv},
	author = {Gabbard, James and van Rees, Wim M.},
	month = sep,
	year = {2023},
}

@article{gabbard_high-order_2024,
	title = {A high-order finite difference method for moving immersed domain boundaries and material interfaces},
	volume = {507},
	issn = {0021-9991},
	url = {https://www.sciencedirect.com/science/article/pii/S0021999124002286},
	doi = {10.1016/j.jcp.2024.112979},
	urldate = {2024-11-19},
	journal = {Journal of Computational Physics},
	author = {Gabbard, James and van Rees, Wim M.},
	month = jun,
	year = {2024},
	pages = {112979},
}

@article{juul_spietz_regularization_2018,
	title = {A regularization method for solving the {Poisson} equation for mixed unbounded-periodic domains},
	volume = {356},
	issn = {0021-9991},
	url = {https://www.sciencedirect.com/science/article/pii/S0021999117309038},
	doi = {10.1016/j.jcp.2017.12.018},
	urldate = {2023-12-12},
	journal = {Journal of Computational Physics},
	author = {Juul Spietz, Henrik and Mølholm Hejlesen, Mads and Walther, Jens Honoré},
	month = mar,
	year = {2018},
	pages = {439--447},
}

@article{deriaz_high-order_2023,
	title = {High-order {Adaptive} {Mesh} {Refinement} multigrid {Poisson} solver in any dimension},
	volume = {480},
	issn = {0021-9991},
	url = {https://www.sciencedirect.com/science/article/pii/S0021999123001079},
	doi = {10.1016/j.jcp.2023.112012},
	urldate = {2025-02-12},
	journal = {Journal of Computational Physics},
	author = {Deriaz, Erwan},
	month = may,
	year = {2023},
	pages = {112012},
}

@article{brown_multigrid_2005,
	title = {Multigrid elliptic equation solver with adaptive mesh refinement},
	volume = {209},
	issn = {0021-9991},
	url = {https://www.sciencedirect.com/science/article/pii/S0021999105001890},
	doi = {10.1016/j.jcp.2005.03.026},
	language = {en},
	number = {2},
	urldate = {2022-10-26},
	journal = {Journal of Computational Physics},
	author = {Brown, J. David and Lowe, Lisa L.},
	month = nov,
	year = {2005},
	pages = {582--598},
}

@article{balty_2025,
    author = {Balty, Pierre and Duponcheel, Matthieu and Chatelain, Philippe},
    title = {An adaptive multiresolution Vortex Particle-Mesh method for the simulation of unbounded incompressible flows},
    journal = {Computer Methods in Applied Mechanics and Engineering},
    year = {accepted},
}

@misc{flups,
    key = {flups},
    title = {{FLUPS} {G}it{H}ub repository},
    howpublished = {\url{https://github.com/vortexlab-uclouvain/flups}},
}

@misc{murphy,
    key = {murphy},
    title = {Murphy {G}it{H}ub repository},
    howpublished = {\url{https://github.com/vanreeslab/murphy}},
}

@article{gholami_fft_2016,
	title = {{FFT}, {FMM}, or {Multigrid}? {A} comparative {Study} of {State}-{Of}-the-{Art} {Poisson} {Solvers} for {Uniform} and {Nonuniform} {Grids} in the {Unit} {Cube}},
	volume = {38},
	issn = {1064-8275},
	shorttitle = {{FFT}, {FMM}, or {Multigrid}?},
	url = {https://epubs.siam.org/doi/abs/10.1137/15M1010798},
	doi = {10.1137/15M1010798},
	number = {3},
	urldate = {2025-01-09},
	journal = {SIAM Journal on Scientific Computing},
	author = {Gholami, Amir and Malhotra, Dhairya and Sundar, Hari and Biros, George},
	month = jan,
	year = {2016},
	pages = {C280--C306},
}

@article{ibeid_fft_2020,
	title = {{FFT}, {FMM}, and multigrid on the road to exascale: {Performance} challenges and opportunities},
	volume = {136},
	issn = {0743-7315},
	shorttitle = {{FFT}, {FMM}, and multigrid on the road to exascale},
	url = {https://www.sciencedirect.com/science/article/pii/S0743731518305513},
	doi = {10.1016/j.jpdc.2019.09.014},
	language = {en},
	urldate = {2022-01-19},
	journal = {Journal of Parallel and Distributed Computing},
	author = {Ibeid, Huda and Olson, Luke and Gropp, William},
	month = feb,
	year = {2020},
	pages = {63--74},
}

@article{bastian1998load,
  title={Load balancing for adaptive multigrid methods},
  author={Bastian, Peter},
  journal={SIAM Journal on Scientific Computing},
  volume={19},
  number={4},
  pages={1303--1321},
  year={1998},
  publisher={SIAM}
}

@article{hou_adaptive_2024,
	title = {An adaptive lattice {Green}'s function method for external flows with two unbounded and one homogeneous directions},
	volume = {519},
	issn = {0021-9991},
	url = {https://www.sciencedirect.com/science/article/pii/S0021999124006181},
	doi = {10.1016/j.jcp.2024.113370},
	urldate = {2025-03-05},
	journal = {Journal of Computational Physics},
	author = {Hou, Wei and Colonius, Tim},
	month = dec,
	year = {2024},
	pages = {113370},
}

@article{liska2014,
title = {A parallel fast multipole method for elliptic difference equations},
journal = {Journal of Computational Physics},
volume = {278},
pages = {76-91},
year = {2014},
issn = {0021-9991},
doi = {https://doi.org/10.1016/j.jcp.2014.07.048},
url = {https://www.sciencedirect.com/science/article/pii/S0021999114005415},
author = {Sebastian Liska and Tim Colonius}
}

@article{berger1958use,
  title={The use of discrete Green's functions in the numerical solution of Poisson's equation},
  author={Berger, JM and Lasher, GJ},
  journal={Illinois Journal of Mathematics},
  volume={2},
  number={4A},
  pages={593--607},
  year={1958},
  publisher={Duke University Press}
}

@article{gillman_fast_2010,
	title = {Fast and accurate numerical methods for solving elliptic difference equations defined on lattices},
	volume = {229},
	issn = {0021-9991},
	url = {https://www.sciencedirect.com/science/article/pii/S0021999110004171},
	doi = {10.1016/j.jcp.2010.07.024},
	number = {24},
	urldate = {2025-03-05},
	journal = {Journal of Computational Physics},
	author = {Gillman, A. and Martinsson, P. G.},
	month = dec,
	year = {2010},
	pages = {9026--9041},
}

@misc{gholami_accfft_2016,
	title = {{AccFFT}: {A} library for distributed-memory {FFT} on {CPU} and {GPU} architectures},
	shorttitle = {{AccFFT}},
	url = {http://arxiv.org/abs/1506.07933},
	language = {en},
	urldate = {2022-07-11},
	publisher = {arXiv},
	author = {Gholami, Amir and Hill, Judith and Malhotra, Dhairya and Biros, George},
	month = may,
	year = {2016},
	note = {arXiv:1506.07933 [cs]},
}

@incollection{krzhizhanovskaya_heffte_2020,
	address = {Cham},
	title = {{heFFTe}: {Highly} {Efficient} {FFT} for {Exascale}},
	volume = {12137},
	isbn = {978-3-030-50370-3 978-3-030-50371-0},
	shorttitle = {{heFFTe}},
	url = {http://link.springer.com/10.1007/978-3-030-50371-0_19},
	language = {en},
	urldate = {2022-03-03},
	booktitle = {Computational {Science} – {ICCS} 2020},
	publisher = {Springer International Publishing},
	author = {Ayala, Alan and Tomov, Stanimire and Haidar, Azzam and Dongarra, Jack},
	editor = {Krzhizhanovskaya, Valeria V. and Závodszky, Gábor and Lees, Michael H. and Dongarra, Jack J. and Sloot, Peter M. A. and Brissos, Sérgio and Teixeira, João},
	year = {2020},
	doi = {10.1007/978-3-030-50371-0_19},
	note = {Series Title: Lecture Notes in Computer Science},
	pages = {262--275},
}

@article{pekurovsky_p3dfft_2012,
	title = {{P3DFFT}: {A} {Framework} for {Parallel} {Computations} of {Fourier} {Transforms} in {Three} {Dimensions}},
	volume = {34},
	issn = {1064-8275},
	shorttitle = {{P3DFFT}},
	url = {https://epubs.siam.org/doi/abs/10.1137/11082748X},
	doi = {10.1137/11082748X},
	number = {4},
	urldate = {2025-03-05},
	journal = {SIAM Journal on Scientific Computing},
	author = {Pekurovsky, Dmitry},
	month = jan,
	year = {2012},
	note = {Publisher: Society for Industrial and Applied Mathematics},
	pages = {C192--C209},
}

@article{pozrikidis2001note,
  title={A note on the regularization of the discrete Poisson--Neumann problem},
  author={Pozrikidis, C},
  journal={Journal of Computational Physics},
  volume={172},
  number={2},
  pages={917--923},
  year={2001},
  publisher={Elsevier}
}

\end{document}